\documentclass[11pt]{article}
\usepackage{amsfonts,mathrsfs}
\usepackage{amssymb}
\usepackage{amsmath}
\usepackage{color}
=0pt
\def\sqr#1#2{{\vcenter{\vbox{\hrule height.#2pt
              \hbox{\vrule width.#2pt height#1pt \kern#1pt \vrule
width.#2pt}
              \hrule height.#2pt}}}}
\def\signed #1{{\unskip\nobreak\hfil\penalty50
              \hskip2em\hbox{}\nobreak\hfil#1
              \parfillskip=0pt \finalhyphendemerits=0 \par}}
\def\endpf{\signed {$\sqr69$}}
\def\dbF{{\mathbb{F}}}

\def\dbH{{\mathbb{H}}}

\def\dbN{{\mathbb{N}}}

\def\dbP{{\mathbb{P}}}

\def\dbR{{\mathbb{R}}}

\def\d{\delta}
\def\e{\varepsilon}

\def\f{\varphi}

\def\om{\omega}

\def\3n{\negthinspace \negthinspace \negthinspace }
\def\2n{\negthinspace \negthinspace }
\def\1n{\negthinspace }
\def\ns{\noalign{\smallskip} }

\def\ds{\displaystyle}
\def\D{\Delta}

\def\Om{\Omega}
\def\om{\omega}
\def\cC{{\cal C}}

\def\cF{{\cal F}}

\def\cL{{\cal L}}

\def\cP{{\cal P}}

\def\cU{{\cal U}}

\def\mE{{\mathbb{E}}}

\def\no{\noindent}

\def\ss{\smallskip}

\def\bs{\bigskip}
\def\q{\quad}
\def\qq{\qquad}

\def\liminf{\mathop{\underline{\rm lim}}}

\def\lan{\big\langle}
\def\ran{\big\rangle}

\def\pa{\partial}

\def\wt{\widetilde}
\def\cd{\cdot}

\def\dim{\hbox{\rm dim$\,$}}

\def\ae{\hbox{\rm a.e.{ }}}

\def\deq{\mathop{\buildrel\D\over=}}

\def\({\Big (}
\def\){\Big )}
\def\[{\Big[}
\def\]{\Big]}

\def\={\buildrel \triangle \over =}

\def\be{\begin{equation}}
\def\bel{\begin{equation}\label}
\def\ee{\end{equation}}
\def\bea{\begin{eqnarray}}
\def\eea{\end{eqnarray}}
\def\bt{\begin{theorem}}
\def\et{\end{theorem}}
\def\bc{\begin{corollary}}
\def\ec{\end{corollary}}
\def\bl{\begin{lemma}}
\def\el{\end{lemma}}
\def\bp{\begin{proposition}}
\def\ep{\end{proposition}}
\def\br{\begin{remark}}
\def\er{\end{remark}}
\def\ba{\begin{array}}
\def\ea{\end{array}}
\def\bde{\begin{definition}}
\def\ede{\end{definition}}

\newtheorem{lemma}{Lemma}[section]
\newtheorem{remark}{Remark}[section]

\newtheorem{theorem}{Theorem}[section]
\newtheorem{corollary}{Corollary}[section]

\newtheorem{definition}{Definition}[section]
\newtheorem{proposition}{Proposition}[section]

\DeclareMathOperator*\lowlim{\underline{lim}}
\DeclareMathOperator*\uplim{\overline{lim}}

\allowdisplaybreaks
 \makeatletter
   
   \@addtoreset{equation}{section}
\makeatother

\begin{document}

\title{\bf Relationships Between the Maximum Principle and Dynamic Programming for Infinite Dimensional Stochastic Control Systems
}

\author{Liangying Chen\footnote{School
of Mathematics, Sichuan University, Chengdu, P. R. China, and Sorbonne Universit\'es, UPMC Univ Paris 06, Paris, France. Email: chenli@ljll.math.upmc.fr. Liangying Chen is supported by the European Union's Horizon 2020 research and innovation programme under the Marie Sklodowska-Curie grant agreement No 945322.}
~~~ and ~~~ Qi L\"{u}\footnote{School of Mathematics, Sichuan
University, Chengdu, P. R. China. Email: lu@scu.edu.cn. Qi L\"u is supported by the NSF of China under
grants 12025105 and 11971334.}}

\date{}

\maketitle

\begin{abstract}

Pontryagin type  maximum principle and Bellman's dynamic programming
principle serve as two of the most important tools in solving optimal control
problems. There is a huge literature on the study of relationship between them. The main purpose of this paper is to investigate the relationships
between Pontryagin type maximum principle and dynamic programming
principle for  control systems governed by stochastic evolution equations in infinite dimensional space,
with the control variables appearing into both the drift and the
diffusion terms. To do so, we first establish dynamic programming
principle for those systems without employing the martingale solutions. Then we establish the desired
relationships in both cases that value function
associated is smooth and nonsmooth. For the nonsmooth case, in
particular, by employing the relaxed transposition solution, we
discover the connection between the superdifferentials and
subdifferentials of value function and the first-order and
second-order adjoint equations.
\end{abstract}

\bs

\no{\bf 2010 Mathematics Subject
Classification}. 93E20.

\bs

\no{\bf Key Words:} Pontryagin type maximum principle, dynamic programming,
stochastic optimal control, stochastic distributed parameter systems

\section{Introduction}\label{intro}

Pontryagin type  maximum principle (PMP, for
short) and Bellman's dynamic programming (DPP, for short) serve as two of the most important tools in solving optimal control
problems. Consequently, it is  natural to ask the followin question:
\ss

{\bf Problem (R)} Are there any relations
between PMP and DPP?

\ss

{\bf Problem (R)} is first studied for control systems govenred by ordinary differential equations in the early 1960s with the assumption that the value function is continuous differentiable (e.g.,\cite{Pontryagin1962}). For a long time, the result was just formal (except for some special cases)
due to the fact that the value function is not smooth enough. In the late 1980s, with the tools of viscosity solution and nonsmooth analysis, results for {\bf Problem (R)} without assuming the smoothness of the value function are establsihed in \cite{Barron1986,Clarke1987,Zhou1990}. Soon after, these results are generalized to systems by partial differential equations and stochastic differential equations (e.g.,\cite{Cannarsa1992,Cannarsa1996,Zhou1991}).  After that, {\bf Problem (R)} is studied extensively for  different control systems. It is impossible to list the related references comprehensively since there are too many. Hence, we refer the readers to \cite{Cannarsa2015,Cernea2005,Frankowska2013,Li1995,Nie2017,Pakniyat2017,Yong1999} and the rich references therein for the study of {\bf Problem (R)}.   
On the other hand, as far as we know, there is no published results concerning {\bf Problem (R)} for  systems governed by stochastic partial differential equations, or more general, stochastic evolution equations (SEEs for short) in infinite dimensional space.

This paper is a first attempt to study {\bf Problem (R)} for controlled SEEs  in infinite dimensional space. More precisely, in this paper, for the first time, we obtain DPP for infinite
dimensional stochastic control systems with deterministic
coefficients under the strong formulation, and then establish
relationships between PMP (in terms of relaxed transposition
solutions) and DPP for such systems.

To be more specific, let us first introduce the control system studied in this paper.  Let $H$ and $\wt H$ be two
separable Hilbert spaces, $T>0$, $(\Omega, \mathcal{F}, \{\cF_t\}_{t\geq 0},
\mathbb{P})$ be a complete filtered probability space, on which an $\wt H$-valued cylindrical
Brownian motion $W(\cdot)$ is defined and $\textbf{F}\deq \{\cF_t\}_{t\geq 0}$ is the natural filtration generated by $W(\cdot)$.    Denote by $\mathbb{F}$ the
progressive $\sigma$-algebra w.r.t. $\textbf{F}$.
Write $ \mathcal{L}_2^0$ for the space of all
Hilbert-Schmidt operators from $\wt H$ to $H$, which is also a separable Hilbert space. Let $A:D(A)\subset H\to H$ be a linear
operator, which generates a $C_{0}$-semigroup $\{S(t)\}_{t\ge 0}$ on $H$.  Let $U$ be a separable metric space with a metric
${\bf d}(\cdot,\cdot)$. 
For $t\in [0,T)$, put
\begin{equation*}
   \mathcal{U}[t,T]\deq\big\{u:[t,T]\times\Omega\to U\big| u \mbox{ is $\textbf{F}$-adapted}\big\}.
\end{equation*}
In this paper, $C$ is a generic constant which may vary from line to line.

The control system is
\begin{equation}\label{system1}
\begin{cases}\ds
dX(t)=(AX(t)+a(t,X(t),u(t)))dt+b(t,X(t),u(t))dW(t), &t \in (0,T],\\
\ns\ds X(0)=\eta \in H,
\end{cases}
\end{equation}
and the cost functional is
\begin{equation}\label{cost1}
\mathcal{J}(\eta;u(\cdot))=\mathbb{E}\Big(\int_0^T
f(t,X(t),u(t))dt+h(X(t))\Big).
\end{equation}

We make the following assumptions for the control system 
\eqref{system1} and cost functional \eqref{cost1}. 

\ss

{\bf (S1)} {\it  Suppose that: i) $a(\cdot, \cdot,\cdot): [0,T]\times H\times U \to H$ is
	$\mathcal{B}([0,T])\otimes \mathcal{B}(H)\otimes
	\mathcal{B}(U)/\mathcal{B}(H)$-measurable and
	$b(\cdot,\cdot,\cdot): [0,T]\times H\times U \to
	\mathcal{L}_2^0$ is $\mathcal{B}([0,T])\otimes \mathcal{B}(H)\otimes
	\mathcal{B}(U)/\mathcal{B}(\mathcal{L}_2^0)$-measurable; ii) for any
	$(t,\eta)\in [0,T]\times H$, the maps $a(t,\eta,\cdot):U\to H $ and
	$b(t,\eta,\cdot):U\to \mathcal{L}_2^0$ are continuous; and iii)
	for any $(t,\eta_1,\eta_2,u)\in [0,T]\times H \times H\times U$,
	\begin{equation*}
		\begin{cases}
			|a(t,\eta_1,u)-a(t,\eta_2,u)|_H\leq C|\eta_1-\eta_2|_H,\\
			|b(t,\eta_1,u)-b (t,\eta_2,u)|_{\mathcal{L}_2^0}\leq C|\eta_1-\eta_2|_H,\\
			|a(t,0,u)|_H \leq C, \q \ \ \ \ |b(t,0,u)|_{\mathcal{L}_2^0}\leq
			C.
		\end{cases}
\end{equation*}}

\ss

{\bf (S2)} {\it  Suppose that: i) $f(\cdot,\cdot,\cdot):[0,T]\times H\times U\to \mathbb{R}$
	is $\mathcal{B}([0,T])\otimes \mathcal{B}(H)\otimes
	\mathcal{B}(U)/\mathcal{B}(\mathbb{R})$-measurable and
	$h(\cdot):H\to \mathbb{R}$ is
	$\mathcal{B}(H)/\mathcal{B}(\mathbb{R})$-measurable; ii) For any
	$(t,\eta)\in [0,T]\times H$, the functional $f(t,\eta,\cdot): U\to
	\mathbb{R}$ is continuous; and iii) For any $(t,\eta_1,\eta_2,u)\in
	[0,T]\times H\times H\times U$,
	\begin{equation*}
		\begin{cases}
			|f(t,\eta_1,u)-f(t,\eta_2,u)|\leq C|\eta_1-\eta_2|_H, \\
			|h(\eta_1)-h(\eta_2)|\leq C|\eta_1-\eta_2|_H\\
			|f(t,0,u)|\leq C,  \q \ \ |h(0)|\leq C.
		\end{cases}
\end{equation*}}

\ss

{\bf (S3)} {\it  The maps $a(t,\eta,u)$ and $b (t,\eta,u)$, and the functionals $f(t,\eta,u)$ and $h(x)$ are $C^2$ with respect to $x$, such that for $\phi(t,\eta,u)=a(t,\eta,u),\ b(t,\eta,u)$ and $\Psi (t,\eta,u)=f(t,\eta,u),\ h(x)$, it holds that $\phi_x(t,\eta,u),\ \Psi_x(t,\eta,u),\ \phi_{xx}(t,\eta,u),$ and $\Psi_{xx}(t,\eta,u)$ are continuous with respect to $u$.  Moreover, there is a  modulus of continuity $\bar{\omega}:[0,\infty)\to [0,\infty)$ such that for any $(t,\eta,,\eta_1,\eta_2,u)\in [0,T]\times H\times H\times H\times U$,
	\begin{equation*}
		\begin{cases}
			|a_{xx}(t,\eta,u)|_{\mathcal{L}(H,H;H)}+|b_{xx}(t,\eta,u)|_{\mathcal{L}(H,H;\mathcal{L}_2^0)}+|\Psi_{xx}(t,\eta,u)|_{\mathcal{L}(H)}\leq C;\\
			|a_{xx}(t,\eta_1,u_1)-a_{xx}(t,\eta_2,u_2)|_{\mathcal{L}(H,H;H)}+|b _{xx}(t,\eta_1,u_1)-b_{xx}(t,\eta_2,u_2)|_ {\mathcal{L}(H,H;\mathcal{L}_2^0)}\\
			\  + |\Psi_{xx}(t,\eta_1,u_1)-\Psi_{xx}(t,\eta_2,u_2)|_{\mathcal{L}(H)} \leq \bar{\omega}(|\eta_1-\eta_2|)+\mathbf d(u_1,u_2).
		\end{cases}
\end{equation*}}

Under {\bf (S1)}, for any $u(\cdot)\in \mathcal{U}[0,T]$,   the control system
\eqref{system1} has a unique mild solution $X(\cdot) \in C_\dbF([0,T];$ $L^2(\Omega;H))$ (see \cite[Theorem 3.14]{Lu2021} for example).

Consider the following optimal control problem:

\textbf{Problem} $\boldsymbol{(S_{\eta})}$.  For any given $\eta\in
H$, find a
$\bar{u}(\cdot)\in \mathcal{U}[0,T]$ such that
\begin{equation}\label{OP1}
\mathcal{J}(\bar{u}(\cdot))=\inf\limits_{u(\cdot)\in \mathcal{U}[0,T]}\mathcal{J}(u(\cdot)). 
\end{equation}

Any $\bar{u}(\cdot)\in \mathcal{U}[0,T]$ satisfying \eqref{OP1} is called an {\it optimal control} (of 
\textbf{Problem} $\boldsymbol{(S_{\eta})}$). The corresponding state
$\overline{X}(\cdot)$ is called an {\it optimal state}, and
$(\overline{X}(\cdot),\bar{u}(\cdot))$ is called an {\it optimal pair}.

To study Problem $\boldsymbol{(S_{\eta})}$, people introduce two tools, that is, PMP (see Section \ref{sec-PMP} for details) and DPP (see Section \ref{sec-DPP} for details). In this paper, we will investigate the relationship between these two tools.  The rest of this paper is as follows. In Section \ref{sec-PMP}, we recall the
PMP for Problem $\boldsymbol{(S_{\eta})}$. In Section \ref{sec-DPP}, unlike the formulation of stochastic dynamic programming principle in the literature, which employs martingale solutions to SEEs,  we establish a
stochastic dynamic programming for the classical mild solution solutions to SEEs.   Section \ref{sec-smooth} is addressed to the relationship between PMP and DPP when the value function is smooth
enough, while Section 
\ref{sec-nonsmooth} is devoted to the same problem for general value functions.

\section{Pontryagin type maximum principle for \textbf{Problem} $\boldsymbol{(S_{\eta})}$ }\label{sec-PMP}

For the convenience of readers, let us recall the PMP for \textbf{Problem} $\boldsymbol{(S_{\eta})}$, which were established recently in several papers under different assumptions (e.g.,\cite{Du2013,Fuhrman2013,Lu2014,Lu2015}). 

To establish the PMP, we need to introduce the first-order
adjoint equation:
\begin{equation}\label{ad-eq1}
\begin{cases}
   dp(t)=-A^*p(t)dt-\big(a_x(t,\overline{X}(t),\bar{u}(t))^*p(t)+b_x(t,\overline{X}(t),\bar{u}(t))^*q(t)\\
   \ \ \ \ \ \ \ \ \q  -f_x(t,\overline{X}(t),\bar{u}(t))\big)dt+q(t)dW(t), \qq \q \q \qq cccc[0,T),\\
   p(T)=-h_x(\overline{X}(T)),
\end{cases}
\end{equation}
and the second-order adjoint equation:
\begin{eqnarray}\label{ad-eq2}
\begin{cases}
   dP(t)=-\big[ \big(A^*+a_x(t,\overline{X}(t),\bar{u}(t))^*\big)P(t)+ P(t)\big(A+a_x(t,\overline{X}(t),\bar{u}(t))\big) \\
   \ \ \ \ \ \ \ \ \ \q +b_x(t,\overline{X}(t),\bar{u}(t))^*P(t)b _x(t,\overline{X}(t),\bar{u}(t))+b_x(t,\overline{X}(t),\bar{u}(t))^*Q(t) \\
   \ \ \ \ \ \ \ \ \ \q  +Q(t)b_x(t,\overline{X}(t),\bar{u}(t))\!+\!\mathbb{H}_{xx}(t,\overline{X}(t),\bar{u}(t),p(t),q(t))\big]dt\! +Q(t)dW(t), \ \  t\in [0,T),\\
   P(T)=-h_{xx}(\overline{X}(T)),
\end{cases}
\end{eqnarray}
with
$$
\begin{array}{ll}\ds
\mathbb{H}(t,\eta,u,p,q)=\langle p,a(t,\eta,u)\rangle_H + \langle q,b(t,\eta,u)\rangle_{\mathcal{L}_2^0}-f(t,\eta,u),\\
\ns\ds \qq\qq\qq\qq (t,\eta,u,p,q)\in[0,T]\times H\times U\times H\times
\mathcal{L}_2^0.
\end{array}
$$

The equation \eqref{ad-eq1} is an $H$-valued backward stochastic evolution equation (BSEE for short). It admits a unique mild solution $(p,q)\in L^2_{\dbF}(\Om;C[0,T];H)\times L^2_{\dbF}(0,T;H)$ (e.g.,\cite[Section 4.2]{Lu2021}). 

When $H = \dbR^n$,  \eqref{ad-eq2} is an $\dbR^{n\times n}$ (matrix)-valued backward stochastic differential equation (which can be easily regarded as an  $\dbR^{n^2}$ (vector)-valued backward stochastic
differential equation), and therefore, its well-posedness follows from
the one for backward stochastic evolution equations valued in Hilbert spaces (e.g., \cite[Section 4.2]{Lu2021}).
When $\dim H =\infty$,  although $L(H)$ is still a
Banach space, it is neither reflexive (needless to say to be a Hilbert space)
nor separable. To the best of our knowledge, in the previous literatures there exists no such a stochastic integration/evolution equation theory
in general Banach spaces that can be employed to treat the well-posedness
of \eqref{ad-eq2} in the usual sense. Then, a new notion of solution, i.e., the relaxed transposition solution, is introduced to \eqref{ad-eq2}.  Let us  briefly recall it. 

For simplicity of notation, we put
$$
\begin{cases}
J(t)=a_x(t,\overline{X}(t),\bar{u}(t)),\q K(t)=b_x(t,\overline{X}(t),\bar{u}(t)), \\
\ns\ds F(t)= -\mathbb{H}_{xx}(t,\overline{X}(t),\bar{u}(t),p(t),q(t)),\q P_T=-h_{xx}(\overline{X}(T)).
\end{cases}
$$
For $1\leq p_1, p_2,q_1,q_2\leq\infty$, let
\begin{equation*} 
	\begin{array}{ll}\ds
		\mathcal{L}_{pd}(L_{\mathbb{F}}^2(0,T;L^{4}(\Omega,H));L^2_{\mathbb{F}}(0,T;L^{\frac{4}{3}}(\Omega,H)))\\
		\ns\ds \deq \Big\{L\!\in\! \cL\big(L_{\mathbb{F}}^2(0,T;L^{4}(\Omega,H));L^2_{\mathbb{F}}(0,T;L^{\frac{4}{3}}(\Omega,H))\big) \big| \mbox{for }\ae (t,\omega)\in [0,T]\times\Omega, \mbox{there
			is }\\
		\ns\ds\q \wt L(t,\omega)\!\in\!\cL (H_1;\!H_2)\;  \mbox{such that } \big( L u(\cd)\big)(t,\omega) 
		= \wt L (t,\omega)v(t,\omega),
		\forall\; v(\cd)\in
		L_{\mathbb{F}}^2(0,T;L^{4}(\Omega,H))\Big\}.
	\end{array}
\end{equation*}
In the sequel, if there is no confusion,  we identify $L\in \mathcal{L}_{pd}(L_{\mathbb{F}}^2(0,T;L^{4}(\Omega,H));L^2_{\mathbb{F}}(0,T; L^{\frac{4}{3}}(\Omega,H)))$ with $\wt L(\cd,\cd)$.

Let
\begin{eqnarray*}
	\cP[0,T] \3n& \deq\3n & \big\{P(\cdot,\cdot)\ |\ P(\cdot,\cdot)\in \mathcal{L}_{pd}(L_{\mathbb{F}}^2(0,T;L^{4}(\Omega,H));L^2_{\mathbb{F}}(0,T;L^{\frac{4}{3}}(\Omega,H))),\\
	& & \ \ \ P(\cdot,\cdot)\xi \in D_{\mathbb{F}}([t,T];L^{\frac{4}{3}}(\Omega,H))) \ \textup{and} \ |P(\cdot,\cdot)\xi|_{D_{\mathbb{F}}([t,T];L^{\frac{4}{3}}(\Omega,H))}\\
	& & \ \ \ \leq C |\xi|_{L_{\mathcal{F}_t}^{4}(\Omega;H)} \ \textup{for every}\ t\in [0,T] \ \textup{and} \ \xi \in L_{\mathcal{F}_t}^{4}(\Omega;H)\big\},
\end{eqnarray*}
and
\begin{eqnarray*}
 \mathcal{Q} [0,T] 
\3n& \deq\3n &\big\{(Q^{(\cdot)},\widehat{Q}^{(\cdot)})\ |\ Q^{(t)},\widehat{Q}^{(t)}\in \mathcal{L}(\mathcal{H}_t;L_{\mathbb{F}}^2(t,T;L^{\frac{4}{3}}(\Omega;\mathcal{L}_2^0)))\\
	& & \hspace{0.8cm} \textup{and}\  Q^{(t)}(0,0,\cdot)^{*}=\widehat{Q}^{(t)}(0,0,\cdot) \ \textup{for any}\  t\in [0,T)\big\} 
\end{eqnarray*}
with
\begin{equation*}
	\mathcal{H}_t\deq L_{\mathcal{F}_t}^{4}(\Omega;H)\times L_{\mathbb{F}}^2(t,T;L^{4}(\Omega;H))\times L_{\mathbb{F}}^2(t,T;L^{4}(\Omega;\mathcal{L}_2^0)),\ \ \forall \  t\in [0,T).
\end{equation*}
For $j=1, 2$ and $t\in [0,T)$, consider the following equation:
\begin{equation}\label{test-eq1}
	\begin{cases}
		d\f_j=(A+J)\f_jds+u_jds+Kx_jdW(s)+v_jdW(s) &\textup{ in } (t,T],\\
		\f_j(t)=\xi_j
	\end{cases}
\end{equation}
where $\xi_j \in L_{\mathcal{F}_t}^{4}(\Omega;H)$, $u_j\in L_{\mathbb{F}}^2(t,T;L^{4}(\Omega;H)) $ and $v_j \in
L_{\mathbb{F}}^2(t,T; L^{4}(\Omega;\mathcal{L}_2^0))$.  By the classical well-posedness result for SEEs, we know that \eqref{test-eq1} has a unique mild solution $\f_j\in C_\dbF([t,T];L^4(\Om;H))$ (e.g.,\cite[Section 3.2]{Lu2021}).

\begin{definition}\label{def2.1} A 3-tuple $(P(\cdot), Q^{(\cdot)},\widehat{Q}^{(\cd)})\in
	\mathcal{P}[0,T]\times \mathcal{Q}[0,T]$ is called a
	relaxed transposition solution to the equation \eqref{ad-eq2}  if for any $t\in
	[0,T]$, $\xi_j \in L_{\mathcal{F}_t}^{4}(\Omega; H)$, $ 
	u_j(\cdot)\in L_{\mathbb{F}}^2(t,T;L^{4}(\Omega;H))$
	and $v_j(\cdot)\in
	L_{\mathbb{F}}^2(t,T;L^{4}(\Omega;\mathcal{L}_2^0))$ ($j=1,2$), it holds that
	\begin{eqnarray*}
		& & \mathbb{E}\lan P_T\f_1(T),\f_2(T)\ran_H-\mathbb{E}\int_t^T \lan F(s)\f_1(s),\f_2(s)\ran_H ds\\
		& & =\mathbb{E}\langle P(t)\xi_1,\xi_2\rangle_H+\mathbb{E}\int_t^T \lan P(s)u_1(s),\f_2(s)\ran_Hds+\mathbb{E}\int_t^T \lan P(s)\f_1(s),u_2(s)\ran_Hds\\
		& & \ \ \ +\mathbb{E}\int_t^T \lan P(s)K(s)\f_1(s),v_2(s)\ran_{\mathcal{L}_2^0}ds +\mathbb{E}\int_t^T \lan P(s)v_1(s),K(s)\f_2(s)+v_2(s)\ran_{\mathcal{L}_2^0}ds\\
		& & \ \ \ +\mathbb{E}\int_t^T \lan v_1(s),\widehat{Q}^{(t)}(\xi_2,u_2,v_2)(s)\ran_{\mathcal{L}_2^0}ds +\mathbb{E}\int_t^T \lan  Q^{(t)}(\xi_1,u_1,v_1)(s),v_2(s)\ran_{\mathcal{L}_2^0}ds
	\end{eqnarray*}
\end{definition}

As an immediate corollary of \cite[Theorem 12.9]{Lu2021}, we have
the following well-posedness result for the equation \eqref{ad-eq2}.

\begin{proposition}\label{prp2.1} The equation \eqref{ad-eq2}
	admits a unique relaxed transposition solution $(P(\cdot),
	Q^{(\cdot)}, \widehat{Q}^{(\cdot)})$. Furthermore,
	\begin{eqnarray}
		& & |P|_{\mathcal{L}(L_{\mathbb{F}}^2(0,T;L^{4}(\Omega,H));L^2_{\mathbb{F}}(0,T;L^{\frac{4}{3}}(\Omega,H)))} + \sup\limits_{t\in[0,T)}|(Q^{(t)},\widehat{Q}^{(t)})|_{ \mathcal{L}(\mathcal{H}_t;L_{\mathbb{F}}^2(t,T;L^{\frac{4}{3}}(\Omega;\mathcal{L}_2^0)))^2}\nonumber\\
		& & \leq C\big(|F|_{L_{\mathbb{F}}^1(0,T;L^2(\Omega;L(H)))}+|P_T|_{L_{\mathcal{F}_T}^2(\Omega;L(H))}\big).
	\end{eqnarray}
\end{proposition}
 ($\overline{X} (\cdot), \bar{u}(\cdot), p(\cdot),
q(\cdot), P(\cdot),
Q^{(\cdot)}, \widehat{Q}^{(\cdot)}$) is called an {\it optimal 7-tuple of \textbf{Problem} $\boldsymbol{(S_{\eta})}$}.
Now we can present the PMP for \textbf{Problem} $\boldsymbol{(S_{\eta})}$.
\begin{theorem}\label{maximum p2-1}
	Suppose that  the assumptions
	{\bf (S1)}--{\bf (S3)} hold. Let $(\overline{X} (\cdot), \bar{u}(\cdot), p(\cdot),
	q(\cdot),P(\cdot),
	Q^{(\cdot)},$ $\widehat{Q}^{(\cdot)})$ be an optimal 6-tuple of \textbf{Problem} $\boldsymbol{(S_{\eta})}$.
	Then,   for  a.e. $(t,\omega)\in [0,T]\times
	\Omega$ and for all $\rho \in U$,
	\begin{equation*}\label{MP2-eq1-1}
		\begin{array}{ll}\ds
			\dbH\big(t,\overline X(t),\bar u(t),p(t),q(t)\big) - \dbH\big(t,\bar y(t),\rho,p(t),q(t)\big) \\
			\ns\ds   - \frac{1}{2}\big\langle
			P(t)\big[ b\big(t,\overline X(t),\bar
			u(t)\big)-b\big(t,\overline X(t),\rho\big)
			\big],  b\big(t,\overline X(t),\bar
			u(t)\big)-b\big(t,\overline X(t),\rho\big)
			\big\rangle_{H} \geq 0.
		\end{array}
	\end{equation*}
\end{theorem}

\section{The Dynamic Programming Principle  for  Problem  $\boldsymbol{(S_{\eta})}$}\label{sec-DPP}

In the literature, stochastic dynamic programming principle for  \textbf{Problem} $\boldsymbol{(S_{\eta})}$ in weak formulation is already established. A nice treatise for that is \cite{Fabbri2017}. In that formulation, probability spaces and Brownian motions vary with the controls. In other words, the probability
space and Brownian motion are  part of the control. Usually, optimal control problems for SEEs are formulated in strong formulation, i.e., the probability space and the Brownian motion are fixed. Hence, it is natural to ask whether the stochastic dynamic programming principle holds in strong formulation. This question is first answered in \cite[Chapter 2, Section 5]{Peng1997} when $H=\dbR^n$. In this section,  we generalize the results in \cite{Peng1997} to \textbf{Problem} $\boldsymbol{(S_{\eta})}$.

First, we introduce a family of optimal control problems. 
For any $(t,\eta)\in [0,T)\times H$,   the control  system is
\begin{equation}\label{system2}
	\begin{cases}
\ds		dX(s)=(AX(s)+a(t,X(s),u(s)))dt+b(s,X(s),u(s))dW(s), &s \in (t,T],\\
\ns\ds		X(t)=\eta,
	\end{cases}
\end{equation}
and the cost functional is
\begin{equation}\label{cost2}
	\mathcal{J}(t,\eta;u(\cdot))=\mathbb{E}\Big(\int_t^T
	f(s,X(s),u(s))ds+h(X(T))\Big).
\end{equation}
For any $u(\cdot)\in \mathcal{U}[t,T]$, it follow 
immediately from the classical well-posedness of SEEs (e.g., \cite[Theorem 3.14]{Lu2021}) that  the control system
\eqref{system2} has a unique mild solution $X(\cdot) \in C_\dbF([t,T];L^2(\Omega;H))$. Hence, the cost functional \eqref{cost2} is well-defined. 

Consider the following optimal control problem:

\textbf{Problem} $\boldsymbol{(S_{t\eta})}$.  For any given $(t,\eta)\in
[0,T]\times H$, find a
$\bar{u}(\cdot)\in \mathcal{U}[t,T]$ such that
\begin{equation}\label{OP2}
	\mathcal{J}(t,\eta;\bar{u}(\cdot))=\inf\limits_{u(\cdot)\in \mathcal{U}[t,T]}\mathcal{J}(t,\eta;u(\cdot)). 
\end{equation}

Any $\bar{u}(\cdot)\in \mathcal{U}[t,T]$ satisfying \eqref{OP2} is called an {\it optimal control} (of 
\textbf{Problem} $\boldsymbol{(S_{t\eta})}$). The corresponding state
$\overline{X}(\cdot)$ is called an {\it optimal state}, and
$(\overline{X}(\cdot),\bar{u}(\cdot))$ is called an {\it optimal pair}.

We have the following  Dynamic Programming Principle.

\begin{theorem}\label{th3.2}  For any $(t,\eta)\in [0,T)\times H$,
	\begin{equation}
		V(t,\eta)
		=\inf_{u(\cdot)\in\mathcal{U}[s,T]}\mathbb{E}\Big(\int_t^{\hat
			t}f(s, X(s;t,\eta,u),u(s))ds +V\big(\hat t,X(\hat
		t;t,\eta,u)\big)\Big),\q\forall\,0\le t\le\hat t
		\le T.
	\end{equation}
\end{theorem}

Further, it holds that
\begin{theorem}\label{th3.3}  
	If $(\overline{X}(\cdot),\bar{u}(\cdot))$ is
	optimal for \textbf{Problem} $\boldsymbol{(S_{t\eta})}$), then
	\begin{equation}
		V(r,\overline{X}(r))=\mathbb{E}\Big\{\int_r^T f(\tau,\overline{X}(\tau),\bar{u}(\tau))d\tau +h(\overline{X}(T))\Big|\mathcal{F}_r\Big\},
		\qquad\mathbb{P}\mbox{-a.s.},\ \forall r\in[t,T].
	\end{equation}
\end{theorem}

Further, we have the following regularity properties for $V$.

\begin{proposition}\label{prop3.4-1} For each $t\in [0,T]$, $\eta$ and
	$\eta^{\prime}\in H$, we have
	\begin{equation}\label{prop3.4-1-eq1}
		|V(t,\eta)|\leq C(1+|\eta|_H)
	\end{equation}
	and
	\begin{equation}\label{prop3.4-1-eq2}
		|V(t,\eta)- V(t,\eta^{\prime})| \leq C|\eta-\eta^{\prime}|_H.
	\end{equation}
\end{proposition}

\begin{proposition}\label{prop3.6-1} 
	The function $V(\cd,\eta)$ is continuous.
\end{proposition}

The proof of Theorems \ref{th3.2}--\ref{th3.3} and Propositions \ref{prop3.4-1}--\ref{prop3.6-1} have their own interest. But they are lengthy and technical. We put them in the next section. The readers may skip them when they read this paper at the first time.

Next, we formally derive a partial
differential equation satisfied by the value function $V(\cdot,\cdot)$.

\begin{proposition}\label{prop3.7} 
	Suppose that $V(\cdot,\cdot)\in
C^{1,2}([0,T];H)$, $V_x:[0,T]\times H\to D(A^{*})$, $A^{*}V_x\in
C([0,T]\times H;H)$. Then $V(\cdot,\cdot)$ is a solution to the
following HJB equation:
\begin{equation}\label{prop3.7-eq1}
\begin{cases}
V_t+\langle A^{*}V_x, \eta\rangle_H +  
\inf\limits_{\rho\in U} G(t,\eta,\rho,V_x,V_{xx})=0, & (t,\eta)\in [0,T]\times H,\\
V(t,\eta)=h(\eta), & \eta\in H
\end{cases}
\end{equation}
where
\begin{equation}\label{prop3.7-eq2}
G(t,\eta,\rho,p,P)=\frac{1}{2}\langle
Pb(t,\eta,\rho),b(t,\eta,\rho)\rangle_{\mathcal{L}_2^0} + \langle p,
a(t,\eta,\rho)\rangle_H-f(t,\eta,\rho),
\end{equation}
$$\hspace{5cm}\forall(t,\eta,\rho,p,P)\in[0,T]\times H\times U\times H\times\mathcal{S}(H).$$
\end{proposition}

{\it Proof}. Let $(t,\eta)\in [0,T]\times H$ be given and for any $u\in U$,
taking the constant control $u(\cdot)=\rho$, we have
\begin{eqnarray*}
 0\3n&\leq\3n& \mathbb{E}\Big(\int_t^T f(r,X(r),\rho)dr+V(t+\varepsilon,X(t+\varepsilon))-V(t,X(t))\Big)\\
  & =\3n&  \mathbb{E}\Big[\int_t^T \big(f(r,X(r),\rho)+V_t(r,X(r))+\langle A^{*}V_x(r,X(r)),X(r)\rangle_H + \langle V_x(r,X(r)),a(r,X(r),\rho)\rangle_H \\
  & & \ \ \ \ \ \ \ \ \ \ \ \ \ +\frac{1}{2}\langle V_{xx}(r,X(r))b(r,X(r),\rho),b(r,X(r),\rho)\rangle_{\mathcal{L}_2^0}\big)dr\Big].
\end{eqnarray*}
Taking expectations, dividing both sides by $\varepsilon$ and sending
$\varepsilon\to 0^+$, one can obtain that
$$
V_t+\langle A^{*}V_x, \eta\rangle_H+G(t,\eta,\rho,V_x,V_{xx})\ge 0 ,\ \ \ \forall \rho\in U
$$
Thus,
$$
V_t+\langle A^{*}V_x, \eta\rangle_H+\inf\limits_{\rho\in U}G(t,\eta,\rho,V_x,V_{xx})\ge 0.
$$
Next, for any $\varepsilon >0$ and $s>0$, there exists a
$u^{\varepsilon,s}(\cdot)\in \mathcal{U}[t,T]$ with the corresponding state $X^{\varepsilon,s}(\cdot)$ such that
\begin{eqnarray*}
   \varepsilon \3n& >\3n& \frac{1}{s} \mathbb{E}\Big(\int_t^{t+s} f(r,X^{\varepsilon,s}(r),u^{\varepsilon,s}(r))dr+V(t+s,X^{\varepsilon,s}(t+s))-V(t,\eta)\Big)\\
  & =\3n& \frac{1}{s} \mathbb{E}\Big[\int_t^{t+s} \Big(f(r,X^{\varepsilon,s}(r),u^{\varepsilon,s}(r))+\langle A^{*}V_x(r,X^{\varepsilon,s}(r)),X^{\varepsilon,s}(r))\rangle_H \\
  & & \q \ \ \ \ \ \ \ \ \ \ \ + \langle V_x(r,X^{\varepsilon,s}(r)),a(r,X^{\varepsilon,s}(r),u^{\varepsilon,s}(r))\rangle_H\\
  & & \q \ \ \ \ \ \ \ \ \ \ \ +\frac{1}{2}\langle V_{xx}(r,X^{\varepsilon,s}(r))b(r,X^{\varepsilon,s}(r),u^{\varepsilon,s}(r)),b(r,X^{\varepsilon,s}(r),u^{\varepsilon,s}(r))\rangle_{\mathcal{L}_2^0}\Big)dr\Big]\\
  &  \ge\3n& \frac{1}{s}\mathbb{E}\[ \int_t^{t+s}\Big(V_t(r,X^{\varepsilon,s}(r))+\langle A^{*}V_x(r,X^{\varepsilon,s}(r)),X^{\varepsilon,s}(r)\rangle_H\\
  & & \q \q \ \ \ \ \ \ \ +\inf\limits_{u\in U}G(r,X^{\varepsilon,s}(r),u,V_x(r,X^{\varepsilon,s}(r)),V_{xx}(r,X^{\varepsilon,s}(r)))\Big)dr\]\\
  & & \to V_t(t,\eta)+\langle A^{*}V_x(t,\eta),\eta)\rangle_H+ \inf\limits_{\rho\in U} G(t,\eta,\rho,V_x(t,\eta),V_{xx}(t,\eta)).
\end{eqnarray*}
 Since
$\varepsilon>0$ is arbitrary, we obtain \eqref{prop3.7-eq1}.
\endpf

\section{A stochastic recursive optimal control problem}

In this section, we consider  a more general optimal control problem, i.e., a stochastic recursive optimal control problem.  Compared with \textbf{Problem} $\boldsymbol{(S_{t\eta})}$, it  is more convenient to be handled since one can employ the theory of BSEEs to study it.  

\subsection{Formulation of the stochastic recursive optimal control problem}
 
For any given  $u\in \mathcal{U}[t,T]$,  
consider the following BSEE:
\begin{eqnarray}\label{bsystem1}
	\begin{cases}
		dY\!(s;t,\zeta,u)
		=\!-g(s,\!X(s;\!t,\!\zeta,\!u),Y(s;\!t,\!\zeta,\!u),Z(s;\!t,\!\zeta,\!u),\!u(s))ds\! +\! Z(s;\!t,\!\zeta,\!u)dW\!(s), \ \,
		s \!\in\! [t,T),\\
		Y(T;t,\zeta,u) = \Phi(X(T;t,\zeta,u)).
	\end{cases}
\end{eqnarray}
Here, $X(\cdot)$ is the mild solution to the equation \eqref{system2} with $X(t)=\zeta\in L^2_{\cF_t}(\Om;H)$, and $\Phi : H\to\mathbb{R}$ and $g : [0,T]\times
H\times\mathbb{R}\times \wt H \times U\to \mathbb{R}$
satisfy the following conditions:

\ss

({\bf S4}) {\it  $g(\cd,\eta,y,z,u)$ is Lebesgue measurable, $g(t,\eta,y,z,\cd)$ is continuous.   For some $L>0$, and all $\eta,\eta^{\prime}\in H$,
$y,y^{\prime}\in\mathbb{R}$, $z,z^{\prime}\in \wt H$,
$u,u^{\prime}\in U$ and a.e. $t\in [0,T]$,
\begin{eqnarray*}
	|g(t,\eta,y,z,u)-g(t,\eta^{\prime},y^{\prime},z^{\prime},u^{\prime})|+|\Phi(\eta)-\Phi(\eta^{\prime})| \leq
	L(|\eta-\eta^{\prime}|_H+|\eta-\eta^{\prime}|+|z-z^{\prime}|_{\wt H}) 
\end{eqnarray*}
and
$$ g(t,\eta,0,0,u) + \Phi(\eta)\leq L(1+|\eta|_H).$$}

By the classical well-posedness result for BSEEs (e.g., \cite[Theorem 4.10]{Lu2021} ), the equation \eqref{bsystem1} admits a unique solution $(Y,Z)\in L^2_\dbF(\Om;C([t,T]))\times L^2_\dbF(t,T;\wt H)$.  

We first give some estimates for solutions to the equations \eqref{system2} and \eqref{bsystem1}.
\begin{proposition}\label{prop3.1} For all $t\in [0,T]$,
	$\zeta,\zeta^{\prime}\in L^2_{\mathcal{F}_t}(\Omega;H)$,
	$u(\cdot),u^{\prime}(\cdot)\in\mathcal{U}[t,T]$,
	\begin{equation}\label{prop3.1-eq1}
		\sup_{t\leq s\leq
			T}\mathbb{E}\big(|X(s)|_H^2\big|\mathcal{F}_t\big)\leq
		C(1+|\zeta|_H^2) 
	\end{equation}
	and
	\begin{eqnarray}\label{prop3.1-eq2}
		\sup\limits_{t\leq s\leq
			T}\mathbb{E}\big(|X(s;t,\zeta,u)\!-\!X(s;t,\zeta',u')|_H^2\big|\mathcal{F}_t\big) \!\leq\! C\[\big|\zeta\!-\!\zeta^{\prime}\big|^2_H\!+\!
		\mathbb{E}\(\int_t^T\!\!\big|u(s)\!-\!u^{\prime}(s)\big|^2_Uds\Big|\mathcal{F}_t\)\],
	\end{eqnarray}
	where the constant $C$ is independent of $t\in [0,T]$.
\end{proposition}

{\it Proof}.  Since
$$
X(s) = S(s-t)\zeta + \int_t^s S(s-\tau)a(\tau,X(\tau),u(\tau))d\tau+ \int_t^s S(s-\tau)b(\tau,X(\tau),u(\tau))dW(\tau),
$$
it follows from ({\bf S1}) that 
$$
\begin{array}{ll}\ds
	\mathbb{E}\big(|X(s)|_H^2\big|\mathcal{F}_t\big)\leq \3n&\ds C\[|\zeta|_H^2 + \mE\(\Big|\int_t^s S(s-\tau)a(\tau,X(\tau),u(\tau))d\tau\Big|_H^2\Big|\cF_t\)\\
	\ns&\ds\q  + \mE\(\Big|\int_t^s S(s-\tau)b(\tau,X(\tau),u(\tau))dW(\tau)\Big|_H^2\Big|\cF_t\)\]\\
	\ns&\ds\leq C\Big\{|\zeta|_H^2 +   \int_t^s \big[ 1+ \mE\big(\big|X(\tau)\big|_H^2\big|\cF_t\big)\big] \Big\}.
\end{array}
$$
This, together with Gronwall's inequality, implies that
$$
\begin{array}{ll}\ds
	\mathbb{E}\big(|X(s)|_H^2\big|\mathcal{F}_t\big)\leq C\big(1+|\zeta|_H^2\big),\q \forall s\in [t,T],
\end{array}
$$
which deduces \eqref{prop3.1-eq1} immediately.

By a similar argument, we have
\begin{eqnarray*}
	& &\3n\3n \sup\limits_{t\leq s\leq
		T}\mathbb{E}\big(|X(s;t,\zeta,u)-X(s;t,\zeta',u')|_H^2\big|\mathcal{F}_t\big)\\
	& &\3n\3n \leq C\mathbb{E}\(|\zeta \!-\! \zeta'|_H^2\! + \! \int_t^T\!\big|a(s,0,u(s))\!-\!a(s,0,u'(s))\big|^2_H ds \! +\!  \int_t^T\!\!\big|b(s,0,u(s))-\!b(s,0, u'(s))\big|^2_{\mathcal{L}_2^0}ds\Big|\mathcal{F}_t\) \\
	& &\3n\3n \leq C\[|\zeta-\zeta^{\prime}|^2_H+
	\mathbb{E}\(\int_t^T|u(s)-u^{\prime}(s)|^2_Uds\Big|\mathcal{F}_t\)\],
\end{eqnarray*}
which concludes \eqref{prop3.1-eq2}.
\endpf

\begin{proposition}\label{prop3.2} For any $\zeta, \zeta' \in L^2_{\mathcal{F}_t}(\Omega;H)$, and $u, 
	u'\in \mathcal{U}[t,T]$,
	\begin{equation}\label{prop3.2-eq1}
		\sup\limits_{t\leq s\leq
			T}\mathbb{E}\(\big|Y(s;t,\zeta,u)\big|^2+\int_t^T\big|Z(s;t,\zeta,u)\big|_{\wt H}^2ds\Big|\mathcal{F}_t\)\leq
		C(1+|\zeta|_H^2),
	\end{equation}
and
	\begin{eqnarray}\label{prop3.2-eq2}
		& & \sup\limits_{t\leq s\leq
			T}\mathbb{E}\(\big|Y(s;t,\zeta,u)-Y(s;t,\zeta^{\prime},u^{\prime})\big|^2+\int_t^T\big|Z(s;t,\zeta,u)-Z(s;t,\zeta^{\prime},u^{\prime})\big|_{\wt H}^2ds\Big|\mathcal{F}_t\)\nonumber\\
		& &  \leq
		C\[|\zeta-\zeta^{\prime}|_H^2+\mathbb{E}\(\int_t^T|u(s)-u^{\prime}(s)|_U^2ds\Big|\mathcal{F}_t\)\].
	\end{eqnarray}
\end{proposition} 

{\it Proof}. By  ({\bf S3}), \eqref{prop3.1-eq1}  and the classical well-posedness result for BSEEs (e.g., \cite[Theorem 4.10]{Lu2021}), 
we have
\begin{eqnarray*}
	& & \sup\limits_{t\leq s\leq
		T}\mathbb{E}\(\big|Y(s;t,\zeta,u)\big|^2+\int_t^T\big|Z(s;t,\zeta,u)\big|_{\wt H}^2ds\Big|\mathcal{F}_t\)\\
	& & \leq C\mathbb{E}\(|\Phi(X(t))|^2 +  \int_t^T \big|g(s,X(s),0,0,u(s))\big|^2ds\Big|\mathcal{F}_t\)\\
	& & \leq C\[1+\sup_{t\leq s\leq
		T}\mathbb{E}\big(|X(s)|_H^2\big|\mathcal{F}_t\big)\]\leq C\big(1+|\zeta|_H^2\big).
\end{eqnarray*}
This gives \eqref{prop3.2-eq1}. Similarly,
\begin{eqnarray*}
	& & \sup\limits_{t\leq s\leq
		T}\mathbb{E}\(\big|Y(s;t,\zeta,u)-Y(s;t,\zeta^{\prime},u^{\prime})|^2+\int_t^T\big|Z(s;t,\zeta,u)-Z(s;t,\zeta^{\prime},u^{\prime})\big|_{\wt H}^2ds\Big|\mathcal{F}_t\) \\
	& & \leq C\mathbb{E}\(\big|\Phi(X(T;t,\zeta,u))-\Phi(X(T;t            ,\zeta^{\prime},u^{\prime}))\big|^2 \\
	& & \qq \ \ +  \int_t^T                       |g(s,X(s;t,\zeta,u),0,0,u(s))-g(s,X(s;t,\zeta',u'),0, 0,u'(s))|^2ds\Big|\mathcal{F}_t\) \\
	& & \leq C\[\sup\limits_{t\leq s\leq
		T}\mathbb{E}\big(\big|X(s;t,\zeta,u)-X(s;t,\zeta',u')\big|_H^2\big|\mathcal{F}_t\big) + \mathbb{E}\(\int_t^T|u(s)-u^{\prime}(s)|^2_Uds\Big|\mathcal{F}_t\)\]\\
	& & \leq C\[|\zeta-\zeta^{\prime}|_H^2+\mathbb{E}\(\int_t^T|u(s)     -u^{\prime}(s)|_U^2ds\Big|\mathcal{F}_t\)\].
\end{eqnarray*}
This deduces \eqref{prop3.2-eq2}.
\endpf

Given the control process $u(\cdot)\in\mathcal{U}[t,T]$, we
introduce the associated cost functional:
\begin{equation}\label{12.13-eq29}
	\wt{\mathcal{J}}(t,\eta;u(\cdot))= Y(t;t,\eta,u),\qquad (t,\eta)\in [0,T]\times H,
\end{equation}
and consider the following stochastic recursive optimal control problem:

\ss
\textbf{Problem} $\boldsymbol{(S_{t\eta})}$.  For any given $(t,\eta)\in
[0,T]\times H$, find a
$\bar{u}(\cdot)\in \mathcal{U}[t,T]$ such that
\begin{equation}\label{OP3}
\wt{\mathcal{J}}(t,\eta;\bar{u}(\cdot))=\operatorname*{ess~inf}\limits_{u(\cdot)\in \mathcal{U}[t,T]}\wt{\mathcal{J}}(t,\eta;u(\cdot)). 
\end{equation}

Since  $\wt{\mathcal{J}}(t,\eta;u(\cdot))$ may be a random variable, so we take $\operatorname*{ess~inf}\limits_{u(\cdot)\in \mathcal{U}[t,T]}$ on the right hand side of \eqref{OP3}. 
At a first glance, one may think that $\operatorname*{ess~inf}\limits_{u(\cdot)\in \mathcal{U}[t,T]}\wt{\mathcal{J}}(t,\eta;u(\cdot))$ is a random variable. Fortunately, we have the following result.
\begin{proposition}\label{prop3.3} 
	Under assumptions {\rm({\bf S1})}  and
	{\rm({\bf S4})}, $\operatorname*{ess~inf}\limits_{u(\cdot)\in \mathcal{U}[t,T]}\wt{\mathcal{J}}(t,\eta;u(\cdot))$ is a
	deterministic function of $(t,\eta)\in
	[0,T]\times H$.
\end{proposition}
Before proving Proposition \ref{prop3.3}, we give some preliminaries.  For each $t>0$,  denote by $\mathbf F^t\deq\{\mathcal{F}^t_s\}_{t\leq s\leq T}$
the natural filtration of the Brownian motion $\{W(s)-W(t)\}_{t\leq
	s\leq T}$. Write $\mathbb{F}^t$ for the
progressive $\sigma$-algebra w.r.t. $\textbf{F}^t$.
Let
$$
\mathcal{U}^t \deq \big\{u(\cdot)\in\mathcal{U}[t,T] \big| \ u(s)\text{
	is $\mathbf F^t$-adapted, } \forall\ t\leq
s\leq T\big\},
$$
and
$$
\mathcal{U}_D^t = \Big\{u(s)=\sum\limits_{j=1}^N
u^j(s)1_{\Omega_j}\Big|\, u^j(s)\in\mathcal{U}^t,
\{\Omega_j\}_{j=1}^N\subset \mathcal{F}_t\ \text {is a partition of }
\Omega\Big\}.
$$
From \cite[Lemma 4.12]{Yong2020}, we know that $\mathcal{U}_D^t$ is dense in $\mathcal{U}[t,T]$.
Next, we give the following Lemma.

\begin{lemma}\label{lm3.1} For any $\eta\in H$, and $\{u^j\}_{j=1}^N\subset \mathcal{U}_D^t$,
	the solutions to \eqref{system2} and \eqref{bsystem1} satisfy
	\begin{equation}\label{lm3.1-eq1}
		\begin{cases}
			X\(\cdot;t,\eta,\sum\limits_{j=1}^N\chi_{\Omega_j}u^j\) =
			\sum\limits_{j=1}^N\chi_{\Omega_j}X(\cdot;t,\eta,u^j),\\
			Y\(\cdot;t,\eta,\sum\limits_{j=1}^N \chi_{\Omega_j}u^j\) = \sum\limits_{j=1}^N\chi_{\Omega_j} Y(\cdot;t,\eta,u^j),\\
			Z\(\cdot;t,\eta,\sum\limits_{j=1}^N\chi_{\Omega_j}u^j\)
			= \sum\limits_{j=1}^N\chi_{\Omega_j}Z(\cdot;t,\eta,u^j).
		\end{cases}
	\end{equation}
\end{lemma}

{\it Proof}. For every $j=1,2,\cdots,N$, we denote
$$
(X^j(s), Y^j(s), Z^j(s)) \equiv (X(s;t,\eta,u^j), Y(s;t,\eta,u^j),
Z(s;t,\eta,u^j)).
$$
Then
\begin{eqnarray}\label{12.13-eq1}
	X^j(s) \3n& =\3n & S(s-t)\eta+\int_t^s S(s-r)a(r,X^j(r),u^j(r))dr \nonumber\\
	& & +\int_t^sS(s-r)b (r,X^j(r),u^j(r))dW(r),\quad
	s\in [t,T],
\end{eqnarray}
and
\begin{equation}\label{12.13-eq2}
	Y^j(s)=\Phi(X^j(T))+\int_s^T g(r,X^j(r),Y^j(r),Z^j(r),u^j(r))dr-\int_s^T Z^j(r)dW(r),\quad
	s\in [t,T].
\end{equation}
Multiply $\chi_{\Omega_j}$ on both sides of \eqref{12.13-eq1} and \eqref{12.13-eq2},
and sum the equations. From the  trivial fact that
$$
\sum\limits_{j=1}^N\chi_{\Omega_j}\varphi(x_j)=\varphi\(\sum_{j=1}^N\chi_{\Omega_j}x_j\), \q \forall \f \in C(H),
$$
we get
\begin{eqnarray*}
	\sum\limits_{j=1}^N \chi_{\Omega_j}X^j(s)
	\3n & = \3n& S(s-t)\eta+\int_t^s
	S(s-r)a\(r,\sum\limits_{j=1}^N\chi_{\Omega_j}X^j(r),\sum\limits_{j=1}^N\chi_{\Omega_j}u^j(r)\)dr\\
	& & +\int_t^s
	S(s-r)b \(r,\sum\limits_{j=1}^N\chi_{\Omega_j}X^j(r),\sum\limits_{j=1}^N\chi_{\Omega_j}u^j(r)\)dW(r),
\end{eqnarray*}
and
\begin{eqnarray*}
	\sum_{j=1}^N\chi_{\Omega_j}Y^j(s) \3n& =\3n & \Phi\(\sum_{j=1}^N\chi_{\Omega_j}X^j(T)\)  - \int_s^T \sum\limits_{j=1}^N\chi_{\Omega_j}Z^j(r)dW(r)\\
	& & + \int_s^T g\(r,\sum\limits_{j=1}^N\chi_{\Omega_j}X^j(r),\sum\limits_{j=1}^N\chi_{\Omega_j}Y^j(r),\sum\limits_{j=1}^N\chi_{\Omega_j}Z^j(r),\sum\limits_{j=1}^N\chi_{\Omega_j}u^j(r)\)dr.
\end{eqnarray*}
Then from the uniqueness of the solution of \eqref{system2} and \eqref{bsystem1}, we get \eqref{lm3.1-eq1}.
\endpf

\ss

{\it Proof of Proposition \ref{prop3.3}}. We divide the proof into three steps.

\textbf{Step\ 1.} In this step, we show that
\begin{equation}\label{12.13-eq4}
	\operatorname*{ess~inf}_{u(\cdot)\in\mathcal{U}[t,T]}\wt {\mathcal{J}}(t,\eta;u(\cdot))=\operatorname*{ess~inf}_{u(\cdot)\in
		\mathcal{U}_D^t}\wt {\mathcal{J}}(t,\eta;u(\cdot)).
\end{equation}
Since $\mathcal{U}_D^t\subset \mathcal{U}[t,T]$, we have
$$
\operatorname*{ess~inf}_{u(\cdot)\in\mathcal{U}[t,T]}\wt {\mathcal{J}}(t,\eta;u(\cdot))\leq
\operatorname*{ess~inf}_{u(\cdot)\in
	\mathcal{U}_D^t}\wt {\mathcal{J}}(t,\eta;u(\cdot)).
$$
Consequently, we only need to prove
\begin{equation}\label{12.13-eq4-1}
	\operatorname*{ess~inf}_{u(\cdot)\in\mathcal{U}_D^t}\wt {\mathcal{J}}(t,\eta;u(\cdot))\leq
	\operatorname*{ess~inf}_{u(\cdot)\in\mathcal{U}[t,T]}\wt {\mathcal{J}}(t,\eta;u(\cdot)).
\end{equation}
For any $\varepsilon>0$,
there exists $\tilde{u}(\cdot)\in \mathcal{U}[t,T]$ such that
$$
\dbP\Big\{\wt {\mathcal{J}}(t,\eta;\tilde{u}(\cdot))<\operatorname*{ess~inf}_{u(\cdot)\in\mathcal{U}[t,T]}
\wt {\mathcal{J}}(t,\eta;u(\cdot))+\varepsilon\Big\}=\delta>0.
$$
From Proposition \ref{prop3.2}, we know that for any $\bar{u}(\cdot)\in
\mathcal{U}_D^t$,
$$
\mathbb{E}\big|Y(t;t,\eta,\bar{u})-Y(t;t,\eta,\tilde{u})\big|^2 \leq
C\mathbb{E}\int_t^T\big|\bar{u}(s)-\tilde{u}(s)\big|_U^2ds.
$$
Since $\mathcal{U}_D^t$ is dense in $\mathcal{U}[t,T]$,  
there exists a sequence
$\{u_n(\cdot)\}_{n=1}^\infty\in\mathcal{U}_D^t$ such that
$$
\lim_{n\to\infty}\mathbb{E} \big|Y(t;t,\eta,u_n)-Y(t;t,\eta,\tilde{u})\big|^2 =0.
$$
Then, there exists a subsequence, without loss of
generality, denoted by $\{u_n(\cdot)\}_{n=1}^\infty$ also, such that
$$
\mathbb{P}\(\bigcap_{m=1}^\infty\bigcup_{N=1}^\infty\bigcap_{n=N}^\infty\Big\{\big|Y(t;t,\eta,u_n)-Y(t;t,\eta,\tilde{u})\big|<\frac{1}{m}\Big\}\)=1,
$$
which implies that
$$\mathbb{P}\(\bigcup_{N=1}^\infty\bigcap_{n=N}^\infty\Big\{\big|Y(t;t,\eta,u_n)-Y(t;t,\eta,\tilde{u})\big|<\frac{1}{m}\Big\}\)=1,
\qquad \forall m\in\mathbb{N}.
$$ 
Therefore,
$$
\lim_{N\to\infty}\mathbb{P}\(\bigcap_{n=N}^\infty\Big\{\big|Y(t;t,\eta,u_n)-Y(t;t,\eta,\tilde{u})\big|<\frac{1}{m}\Big\}\)=1,\qquad \forall m\in\mathbb{N},
$$
which deduces that
$$
\lim_{n\to\infty}\mathbb{P}\Big\{\big|Y(t;t,\eta,u_n)-Y(t;t,\eta,\tilde{u})\big|<\frac{1}{m}\Big\}=1,\qquad
\forall m\in\mathbb{N}.
$$ 
Let us choose $m\in\dbN$ to be large enough such that
$1/m<\varepsilon$ and set
\begin{eqnarray*}
	& & \widetilde{\Omega}  = \big\{\omega\in\Omega \big|
	Y(t;t,\eta,\tilde{u})<\operatorname*{ess~inf}_{u(\cdot)\in\mathcal{U}[t,T]}\mathcal{J}(t,\eta;u(\cdot))+\varepsilon\big\};\\
	& & \Omega_n = \Big\{\omega\in\Omega \big|
	\big|Y(t;t,\eta,u_n)-Y(t;t,\eta,\tilde{u})\big|\leq\frac{1}{m}\Big\},\quad n=1,2,\cdots 
\end{eqnarray*}
Then,  $\mathbb{P}(\widetilde{\Omega})=\delta>0$
and $\lim\limits_{n\to\infty}\mathbb{P}(\Omega_n)=1$. We select $n\in\dbN$
large enough such that $\mathbb{P}(\Omega_n)>1-\delta$, then
$$
\mathbb{P}(\widetilde{\Omega}\cap\Omega_n)=\mathbb{P}(\widetilde{\Omega})+\mathbb{P}(\Omega_n)-\mathbb{P}(\widetilde{\Omega}\cup
\Omega_n)>\delta+(1-\delta)-1=0.
$$ 
It is easy to see that
$$
\mathbb{P}\Big\{Y(t;t,\eta,u_n)
<\operatorname*{ess~inf}_{u(\cdot)\in\mathcal{U}[t,T]}\wt {\mathcal{J}}(t,\eta;u(\cdot))+2\varepsilon\Big\}\geq
\mathbb{P}(\tilde{\Omega}\cap\Omega_n)>0.
$$ 
This implies
$$\operatorname*{ess~inf}_{u(\cdot)\in\mathcal{U}_D^t}\wt {\mathcal{J}}(t,\eta;u(\cdot))\leq
\operatorname*{ess~inf}_{u(\cdot)\in\mathcal{U}[t,T]}\wt {\mathcal{J}}(t,\eta;u(\cdot))+2\varepsilon.
$$
From the arbitrariness of $\varepsilon>0$, we get \eqref{12.13-eq4-1}.

\ss

\textbf{Step\ 2.} In this step, we prove
\begin{equation}\label{12.13-eq5}
	\operatorname*{ess~inf}_{u(\cdot)\in\mathcal{U}_D^t}\wt {\mathcal{J}}(t,\eta;u(\cdot))=\operatorname*{ess~inf}_{u(\cdot)\in\mathcal{U}^t}\wt {\mathcal{J}}(t,\eta;u(\cdot))
\end{equation}
Since $\mathcal{U}^t\subset \mathcal{U}_D^t$, we have 
\begin{equation}\label{12.13-eq5-1}
	\operatorname*{ess~inf}_{u(\cdot)\in\mathcal{U}_D^t}\wt {\mathcal{J}}(t,\eta;u(\cdot))\leq \operatorname*{ess~inf}_{u(\cdot)\in\mathcal{U}^t}\wt {\mathcal{J}}(t,\eta;u(\cdot)).
\end{equation}

Now we show the inverse inequality of \eqref{12.13-eq5-1}.
For all $u(\cdot)\in\mathcal{U}_D^t$, we have
$$
\wt {\mathcal{J}}(t,\eta;u(\cdot))=\wt {\mathcal{J}}\(t,\eta;\sum\limits_{j=1}^N\chi_{\Omega_j}u^j(\cdot)\)=\sum\limits_{j=1}^N\chi_{\Omega_j}\wt {\mathcal{J}}(t,\eta;u^j(\cdot)).
$$
For $j=1,2,\cdots,N$, noting that $u^j(\cdot)$ is $\dbF^t$ measurable, we find that $\wt {\mathcal{J}}(t,\eta;u^j(\cdot))$ 
is deterministic. Without loss of generality, we
assume that
$$
\wt {\mathcal{J}}(t,\eta;u^1(\cdot))\leq \wt {\mathcal{J}}(t,\eta;u^j(\cdot)),\quad
\forall j=2,3,\cdots,N.
$$ 
Thus, it holds that
$$
\wt {\mathcal{J}}(t,\eta;u(\cdot))\geq \wt {\mathcal{J}}(t,\eta;u^1(\cdot))\geq
\operatorname*{ess~inf}_{u(\cdot)\in\mathcal{U}^t}\wt {\mathcal{J}}(t,\eta;u(\cdot)).
$$
From the arbitrariness of $u(\cdot)$, we get
$$
\operatorname*{ess~inf}_{u(\cdot)\in\mathcal{U}_D^t}\wt {\mathcal{J}}(t,\eta;u(\cdot))\geq
\operatorname*{ess~inf}_{u(\cdot)\in\mathcal{U}^t}\wt {\mathcal{J}}(t,\eta;u(\cdot)).
$$

\textbf{Step\ 3.} We finish the proof in this step.

From \eqref{12.13-eq4} and \eqref{12.13-eq5}, we see that
\begin{equation}\label{12.13-eq6}
\operatorname*{ess~inf}\limits_{u(\cdot)\in\mathcal{U} }\wt {\mathcal{J}}(t,\eta;u(\cdot))=\operatorname*{ess~inf}\limits_{u(\cdot)\in\mathcal{U}^t}\wt {\mathcal{J}}(t,\eta;u(\cdot)).
\end{equation}
The right hand side of \eqref{12.13-eq6} is deterministic, so is  the left hand side of \eqref{12.13-eq6} is deterministic. This completes the  proof.
\endpf

By Proposition \ref{prop3.3}, we know that   $\operatorname*{ess~inf}\limits_{u(\cdot)\in \mathcal{U}[t,T]}$ on the right hand side of \eqref{OP3} can be replaced by   $\inf\limits_{u(\cdot)\in \mathcal{U}[t,T]}$. Let us define the value function of the stochastic optimal control
problem as follows:
\begin{equation}\label{12.13-eq30}
	\wt V(t,\eta)=\inf_{u(\cdot)\in \mathcal{U}[t,T]}
	\wt{\mathcal{J}}(t,\eta;u(\cdot)),\qquad (t,\eta)\in [0,T]\times H.
\end{equation}

\subsection{Some properties of $\wt V(\cd,\cd)$}\label{ssec-V}

In this subsection, we  
present some properties of $\wt V(\cd,\cd)$, which will be used to establish the DPP for the stochastic recursive optimal control problem.  
The first result reveals some reglarity of  $\wt V(\cd,\cd)$ with
respect to $\eta$. 

\begin{proposition}\label{prop3.4} Suppose {\rm({\bf S1})}  and
	{\rm({\bf S4})} hold. For each $t\in [0,T]$, $\eta$ and
	$\eta^{\prime}\in H$, we have
	\begin{equation}\label{prop3.4-eq1}
		|\wt V(t,\eta)|\leq C(1+|\eta|_H)
	\end{equation}
	and
	\begin{equation}\label{prop3.4-eq2}
		|\wt V(t,\eta)-\wt V(t,\eta^{\prime})| \leq C|\eta-\eta^{\prime}|_H.
	\end{equation}
\end{proposition}

{\it Proof}. From Propositions \ref{prop3.1} and \ref{prop3.2}, for  
$u(\cdot)\in\mathcal{U}[t,T]$, we have
\begin{equation}\label{12.13-eq7}
	\big|\wt {\mathcal{J}}(t,\eta;u(\cdot))\big|\leq C\big(1+|\eta|_H\big),\q \forall \eta\in H
\end{equation}
and
\begin{equation}\label{12.13-eq8}
	\big|\wt {\mathcal{J}}(t,\eta;u(\cdot))-\wt {\mathcal{J}}(t,\eta';u(\cdot))\big|^2\leq
	C|\eta-\eta'|_H^2,\q \forall \eta,\eta'\in H.
\end{equation}
On the other hand, for each $\varepsilon>0$, there exist $u(\cdot)$
and $u^{\prime}(\cdot)\in\mathcal{U}[t,T]$ such that
\begin{eqnarray*}
	& & \wt {\mathcal{J}}(t,\eta;u^{\prime}(\cdot)) - \varepsilon \leq \wt V(t,\eta)\leq
	\wt {\mathcal{J}}(t,\eta;u(\cdot)),\\
	& & \wt {\mathcal{J}}(t,\eta';u(\cdot)) - \varepsilon\leq \wt V(t,\eta')\leq
	\wt {\mathcal{J}}(t,\eta';u'(\cdot)).
\end{eqnarray*}
Then from \eqref{12.13-eq7}, we get
$$
-C\big(1+|\eta|_H\big) - \varepsilon \leq \wt {\mathcal{J}}(t,\eta;u^{\prime}(\cdot)) -
\varepsilon \leq \wt V(t,\eta)\leq \wt {\mathcal{J}}(t,\eta;u(\cdot))\leq
C\big(1+|\eta|_H\big).
$$
From the arbitrariness of $\varepsilon>0$, we obtain \eqref{prop3.4-eq1}.

Similarly,
$$ 
\wt {\mathcal{J}}(t,\eta;u^{\prime}(\cdot))-\wt {\mathcal{J}}(t,\eta^{\prime};u^{\prime}(\cdot))-\varepsilon\leq
\wt V(t,\eta)-\wt V(t,\eta')\leq \wt {\mathcal{J}}(t,\eta;u(\cdot))
-\wt {\mathcal{J}}(t,\eta^{\prime};u(\cdot))+\varepsilon.
$$ 
Thus,
\begin{eqnarray*}
	\big|\wt V(t,\eta)-\wt V(t,\eta^{\prime})\big| \leq \max\big\{\big|\wt {\mathcal{J}}(t,\eta;u(\cdot))-\wt {\mathcal{J}}(t,\eta^{\prime};u(\cdot))\big|,\,\big|\wt {\mathcal{J}}(t,\eta;u^{\prime}(\cdot))
	-v(t,\eta^{\prime};u^{\prime}(\cdot))\big|\big\}+\varepsilon,
\end{eqnarray*}
which, together with \eqref{12.13-eq8}, implies that
\begin{eqnarray*}
	\big|\wt V(t,\eta)-\wt V(t,\eta^{\prime})\big|^2\3n
	&   \leq \3n& 
	2\max\big\{\big|\wt {\mathcal{J}}(t,\eta;u(\cdot))-\wt {\mathcal{J}}(t,\eta^{\prime};u(\cdot))\big|^2,\big|\wt {\mathcal{J}}(t,\eta;u'(\cdot))-\wt {\mathcal{J}}(t,\eta^{\prime};
	u^{\prime}(\cdot))\big|^2\big\}+2\varepsilon^2\\
	&   \leq \3n& 2C\big|\eta-\eta^{\prime}\big|_H^2 +2\varepsilon^2.
\end{eqnarray*}
From the arbitrariness of $\varepsilon>0$, we obtain \eqref{prop3.4-eq2}.
\endpf

In the following, we will define $\wt V(t,\cd)$ for any
$\zeta$ in $L_{\mathcal{F}_t}^2(\Omega;H)$. Before this, we need
to define $\wt {\mathcal{J}}(t,\zeta;u(\cdot))$ for $\zeta$ in
$L_{\mathcal{F}_t}^2(\Omega;H)$.

First, for each simple function
$\zeta=\sum\limits_{i=1}^N 1_{\Omega_i}\eta_i$, where $\{\eta_i\}_{i=1}^N\subset H$, the similar argument
as Lemma \ref{lm3.1} leads to
\begin{equation}\label{12.13-eq32}
	\wt {\mathcal{J}}(t,\zeta;u(\cdot))=Y(t;t,\zeta;u)=\sum\limits_{i=1}^N 1_{\Omega_i}Y(t;t,\eta_i,u)=\sum\limits_{i=1}^N 1_{\Omega_i}\wt {\mathcal{J}}(t,\eta_i;u(\cdot)).
\end{equation}
Hence, we can put
\begin{equation}\label{12.13-eq32-1}
	\wt V(t,\zeta)\deq \sum\limits_{i=1}^N 1_{\Omega_i}\wt V(t,\eta_i)=\sum\limits_{i=1}^N 1_{\Omega_i}\inf_{u(\cdot)\in\mathcal{U}[t,T]}\wt{\mathcal{J}}(t,\eta_i;u(\cdot)).
\end{equation}

Secondly, for any $\zeta$ in $L_{\mathcal{F}_t}^2(\Omega;H)$, there
exists a sequence of simple functions $\{\zeta_j\}_{j=1}^\infty$ 
such that $\zeta=\lim\limits_{j\to \infty} \zeta_j $ in
$L^2_{\cF_t}(\Omega;H)$. From Proposition \ref{prop3.4},  we see that $\{\wt V(t,\zeta_j)\}_{j=1}^\infty$ is a Cauchy sequence in $L^2_{\cF_t}(\Omega)$. Hence, we can define 
\begin{equation}\label{12.13-eq32-2}
	\wt V(t,\zeta)\deq  \lim_{j\to\infty}\wt V(t,\zeta_j)\q \mbox{ in }L^2_{\cF_t}(\Omega).
\end{equation}
Assume $\{\tilde\zeta_j\}_{j=1}^\infty$ is another sequence of simple functions 
such that $\zeta=\lim\limits_{j\to \infty} \tilde\zeta_j $ in
$L^2_{\cF_t}(\Omega;H)$.  By Proposition \ref{prop3.4}, we know that 
$$
\lim_{j\to\infty}\wt V(t,\zeta_j)=  \lim_{j\to\infty}\wt V(t,\tilde\zeta_j).
$$
Hence, $\wt V(t,\zeta)$ is indepedent of the choice of the sequence of simple functions.

\begin{proposition}\label{prop3.5} Suppose {\rm({\bf S1})}  and
	{\rm({\bf S4})} hold.  Fix $t\in[0,T)$ and $\zeta\in
	L^2_{\mathcal{F}_t}(\Omega;H)$. For each $u(\cdot)\in
	\mathcal{U}[t,T]$,   we have
	\begin{equation}\label{prop3.5-eq1}
		\wt V(t,\zeta)\leq Y(t;t,\zeta,u).
	\end{equation}
	On the other hand, for each $\varepsilon>0$, there exists an
	admissible control $u_\e(\cdot)\in\mathcal{U}[t,T]$ such that
	\begin{equation}\label{prop3.5-eq2}
		\wt V(t,\zeta)\geq Y(t;t,\zeta,u_\e)-\varepsilon, \q \mathbb{P} \mbox{-a.s.}
	\end{equation}
\end{proposition}

{\it  Proof}.  We first prove \eqref{prop3.5-eq1}. When  $\zeta=\sum\limits_{i=1}^N1_{\Om_i}\eta_i$ for $\{\Om_i\}_{i=1}^N\subset \cF_t$ and $\{\eta_i\}_{i=1}^N\subset H$, 
for all $u(\cdot)\in \mathcal{U}[t,T]$, we have
$$
Y(t;t,\zeta,u)=Y\(t;t,\sum\limits_{i=1}^N
1_{\Om_i}\eta_i,u\)=\sum\limits_{i=1}^N1_{\Om_i}Y(t;t,\eta_i,u)\geq
\sum_{i=1}^N 1_{\Om_i}\wt V(t,\eta_i)=\wt V(t,\zeta).
$$
When $\zeta\in L^2_{\mathcal{F}_t}(\Omega;H)$, we can choose a
sequence of simple functions $\{\zeta_j\}_{j=1}^\infty$ 
converging to $\zeta$ in $L^2_{\mathcal{F}_t}(\Omega;H)$.
By Proposition \ref{prop3.2}, we have that
\begin{eqnarray*}
	\lim_{j\to \infty}\mathbb{E}
	\big|Y(t;t,\zeta,u)-Y(t;t,\zeta_j,u)\big|^2 = 0.
\end{eqnarray*}
From the definition of $\wt V(t,\zeta)$, we can find that
\begin{eqnarray*}
	\lim_{j\to \infty}\mathbb{E}\big|\wt V(t,\zeta)-\wt V(t,\zeta_j)\big|^2=0.
\end{eqnarray*}
Then, there exists a subsequence $\{\zeta_{j_k}\}_{k=1}^\infty$ of $\{\zeta_j\}_{j=1}^\infty$,   such that
\begin{eqnarray*}
	\begin{cases}\ds
		\lim_{k\to\infty}Y(t;t,\zeta_{j_k},u)=Y(t;t,\zeta,u),  & \mathbb{P} \mbox{-a.s.,}\\
		\ns\ds\lim_{k\to\infty}\wt V(t,\zeta_{j_k})=\wt V(t,\zeta), & \mathbb{P}\mbox{-a.s.} 
	\end{cases}
\end{eqnarray*}
This, together with  $Y(t;t,\zeta_{j_k},u)\geq \wt V(t,\zeta_{j_k})$, $k=1,2,\cdots$,  implies
\eqref{prop3.5-eq1}.

Now we turn to prove \eqref{prop3.5-eq2}. 
Since $H$ is separable, there exists a dense subset $\{\xi_j\}_{j=1}^{\infty}\subset H$ of $H$. Put
$$
\wt{\Om}_j=\Big\{\omega \in \Omega\Big||\zeta(\omega)-\xi_j|_H<\frac{\varepsilon}{3\sqrt{C}}\Big\},\q \forall\ j \in \mathbb{N},
$$
where $C$ is the larger one of the constants in \eqref{prop3.2-eq2} and  \eqref{prop3.4-eq2}. Let
$$
\Om_1=\wt{\Om}_1,\ \Om_k=\wt{\Om}_k\backslash \bigcup_{j=1}^{k-1}\wt{\Om}_j,\q k\in \mathbb{N}.
$$
Then, $\{\Om_j\}_{j=1}^{\infty}\subset\cF_t$ is a partition of $\Omega$. For $\eta=\sum\limits_{j=1}^{\infty}1_{\Om_j}\xi_j$, from \eqref{prop3.2-eq2} and  \eqref{prop3.4-eq2}, we have
\begin{equation*}
	\big|Y(t;t,\zeta,u)-Y(t;t,\eta,u)\big| \leq
	\frac{\varepsilon}{3},\quad \big|\wt V(t,\zeta)-\wt V(t,\eta)\big| \leq
	\frac{\varepsilon}{3}, \q \dbP\mbox{-a.s}.
\end{equation*}
For each $\xi_j$, by Proposition \ref{prop3.3},  we can choose  $u^j\in \mathcal{U}^t$ such that
$$
\wt V(t,\xi_j)\geq Y(t;t,\xi_j,u^j)-\frac{\varepsilon}{3}.
$$
Let $ u(\cdot)=\sum\limits_{j=1}^{\infty}1_{\Om_j}u^j(\cdot)$. 
Then
\begin{eqnarray*}
	Y(t;t,\zeta,u) \3n& \leq\3n & |Y(t;t,\zeta,u)-Y(t;t,\eta,u)|+Y(t;t,\eta,u)\\
	& \leq\3n &
	\frac{\varepsilon}{3}+\sum\limits_{j=1}^{\infty} 1_{\Om_j}Y(t;t,\xi_j,u^j)\\
	& \leq\3n & \frac{\varepsilon}{3} +\sum_{j=1}^{\infty}
	1_{\Om_j}\(\wt V(t,\xi_j)+\frac{\varepsilon}{3}\)
	= \frac{2}{3} \varepsilon +\wt V(t,\eta)\\
	& \leq\3n & \varepsilon+\wt V(t,\zeta), \q \dbP\mbox{-a.s.},
\end{eqnarray*}
which illustrates \eqref{prop3.5-eq2}.
This completes the proof.
\endpf

Next, we devote ourselves to obtaining the continuity of $\wt V(t,\eta)$
with respect to $t$.

\begin{proposition}\label{prop3.6} 
Suppose {\rm({\bf S1})}  and
{\rm({\bf S4})} hold. 	The function $\wt V(\cd,\eta)$ is continuous.
\end{proposition}

{\it Proof}. We define $Y(s;t,\eta,u)$ for all $s\in [0,T]$ by choosing
$Y(s;t,\eta,u)\equiv Y(t;t,\eta,u)$ for $0\leq s\leq t$.
Fix $\eta\in H$, for $0\leq t_1\leq t_2\leq T$ and every $ \varepsilon>0$,
there exist $u_1(\cdot), u_2(\cdot)\in\mathcal{U}[t,T]$, such that
\begin{eqnarray*}
	& & Y(t_1;t_1,\eta;u_2)-\varepsilon\leq \wt V(t_1,\eta)\leq
	Y(t_1;t_1,\eta,u_1),\\
	& & Y(t_2;t_2,\eta,u_1)-\varepsilon\leq \wt V(t_2,\eta)\leq
	Y(t_2;t_2,\eta,u_2).
\end{eqnarray*}
Then,
\begin{eqnarray*}
	Y(t_1;t_1,\eta,u_2)-Y(t_2;t_2,\eta,u_2)-\varepsilon \leq \wt V(t_1,\eta)-\wt V(t_2,\eta) \leq Y(t_1;t_1,\eta,u_1)-Y(t_2;t_2,\eta,u_1)+
	\varepsilon ,
\end{eqnarray*}
which implies that
\begin{eqnarray*}
	|\wt V(t_1,\eta)\!-\!\wt V(t_2,\eta)| \!\leq\!\max\{|Y(t_1;t_1,\eta,u_1)\!-\!Y(t_2;t_2,\eta,u_1)|,
	|Y(t_1;t_1,\eta,u_2)\!-\!Y(t_2;t_2,\eta,u_2)|\}+\varepsilon.
\end{eqnarray*}
Here we only estimate $|Y(t_1;t_1,\eta,u_1)-Y(t_2;t_2,\eta,u_1)|$ and the
estimate of $|Y(t_1;t_1,\eta,u_2)-Y(t_2;t_2,$ $\eta,u_2)|$ is the same. From
Proposition \ref{prop3.2}, we have
\begin{eqnarray}\label{12.13-eq13}
	& & |Y(t_1;t_1,\eta,u_1)-Y(t_2;t_2,\eta,u_1)|^2\nonumber\\
	& & =|Y(0;t_1,\eta,u_1)-Y(0;t_2,\eta,u_1)|^2\nonumber\\
	& &  \leq C \mathbb{E}\big|\Phi(X(T;t_1,\eta,u_1))-\Phi(X(T;t_2,\eta,u_1))\big|^2 \\
	& & \quad +C \mathbb{E}\Big(\int_0^T|1_{[t_1,
		T]}g(s,X(s;t_1,\eta,u_1),Y(s;t_1,\eta,u_1),Z(s;t_1,\eta,u_1),u_1(s)) \nonumber\\
	& & \ \qquad \qquad -1_{[t_2,
		T]}g(s,X(s;t_2,\eta,u_1),Y(s;t_1,\eta,u_1),Z(s;t_1,\eta,u_1),u_1(s))|ds\Big)^2\nonumber\\
	& & = I_1+I_2.\nonumber
\end{eqnarray}
From Propositions \ref{prop3.1} and \ref{prop3.2}, we
get
\begin{equation}\label{12.13-eq14}
	I_1\leq
	C \mathbb{E}|X(T;t_1,\eta,u_1)-X(T;t_2,\eta,u_1)|_H^2\leq C
	\mathbb{E}|X(t_2;t_1,\eta,u_1)-\eta|_H^2 
\end{equation}
and
\begin{equation}\label{12.13-eq15}
	I_2\leq C(t_2-t_1).
\end{equation}
From \eqref{12.13-eq13}--\eqref{12.13-eq15}, we know that
$$|Y(t_1;t_1,\eta,u_1)-Y(t_2;t_2,\eta,u_1)|\leq C(t_2-t_1)^{1/2}+C\big(\mathbb{E}|X(t_2;t_1,\eta,u_1)-\eta|_H^2\big)^{1/2}.$$
The same argument used to $|Y(t_1;t_1,\eta,u_2)-Y(t_2;t_2,\eta,u_2)|^2$
leads to
$$
\big|\wt V(t_1,\eta)-\wt V(t_2,\eta)\big|\leq C(t_2-t_1)^{1/2}+C(\mathbb{E}|X(t_2;t_1,\eta,u_1)-\eta|_H^2)^{1/2}+\varepsilon.
$$
Due to the arbitrariness of $\varepsilon$, we get
$$
\big|\wt V(t_1,\eta)-\wt V(t_2,\eta)\big|\leq C(t_2-t_1)^{1/2}+C\big(\mathbb{E}|X(t_2;t_1,\eta,u_1)-x|_H^2\big)^{1/2}.
$$
From the continuity of $x(\cdot)$ with respect to $t$, we get the
continuity of $\wt V(\cd,\eta)$. The proof is completed.
\endpf

\subsection{DPP for the stochastic recursive optimal control problem}

Now, borrowing some ideas in \cite{Peng1997}, we study the (generalized) dynamic programming
principle for our recursive optimal control problem \eqref{12.13-eq30}. 

Given the initial condition $(t,\eta)$, an admissible control
$u(\cdot)\in\mathcal{U}[t,T]$, a positive number $\delta\leq T-t$
and a real-valued random variable $\zeta\in
L^2_{\mathcal{F}_{t+\delta}}(\Omega)$, we denote
$$ G^{t,\eta;u}_{t,t+\delta}[\zeta]=Y(t),$$
where $(Y(s),Z(s))_{t\leq s\leq t+\delta}$ is the solution of the
following BSDE with time horizon $t+\delta$
\begin{equation*}
	Y(s) =\zeta+\int_s^{t+\delta}g(r,X(r;t,\eta,u),Y(r),Z(r),u(r))dr
	-\int_s^{t+\delta}Z(r)dW(r),\quad t\leq s\leq t+\delta.
\end{equation*}
\begin{theorem}\label{th3.1} Under the assumptions  {\rm ({\bf S1})} and {\rm ({\bf S4})}, the value function $\wt V(\cd,\cd)$ obeys the following dynamic
	programming principle: For each $0<\delta\leq T-t$,
	\begin{equation}\label{th3.1-eq1}
	\wt V(t,\eta)=\operatorname*{ess~inf}_{u(\cdot)\in\mathcal{U}[t,T]}
		G^{t,\eta;u}_{t,t+\delta}[\wt V(t+\delta,X(t+\delta;t,\eta,u))].
	\end{equation}
\end{theorem}

{\it Proof}. By the uniqueness of the solution to \eqref{bsystem1}, we have
$$
\begin{array}{ll}\ds
	Y(s)\3n&\ds = \Phi(X(T;t,\eta,u))+\int_s^{T}g(r,X(r;t,\eta,u),Y(r),Z(r),u(r))dr
	-\int_s^{T}Z(r)dW(r)\\
	\ns&\ds = Y(t+\delta;t,\eta,u)+\int_s^{t+\delta}g(r,X(r;t,\eta,u),Y(r),Z(r),u(r))dr
	-\int_s^{t+\delta}Z(r)dW(r).
\end{array}
$$
Hence,
\begin{equation}\label{12.13-eq11}
	G^{t,\eta;u}_{t,T}\big[\Phi(X(T;t,\eta,u))\big]=G^{t,\eta;u}_{t,t+\delta}\big[Y(t+\delta;t,\eta,u)\big].
\end{equation}
By the uniqueness of the solution to \eqref{system2} and \eqref{bsystem1}, for $s\geq t+\delta$, we have
$$
\begin{array}{ll}\ds
	Y(s)\!\3n&\ds = \Phi(X(T;t,\eta,u))+\int_s^{T}g(r,X(r;t,\eta,u),Y(r),Z(r),u(r))dr
	-\int_s^{T}Z(r)dW(r)\\
	\ns&\ds = \Phi(X(T;t\!+\!\d,X\!(t\!+\!\d),u))\!+\!\int_s^{T}\!\!g(r,X(r;t,X(t\!+\!\d),u),Y(r),Z(r),u(r))dr \!
	-\!\int_s^{T}\!\!Z(r)dW(r).
\end{array}
$$
Consequently,
\begin{equation}\label{12.13-eq12}
	G^{t,\eta;u}_{t,t+\delta}\big[Y(t+\delta;t,\eta,u)\big]=G^{t,\eta;u}_{t,t+\delta}\big[Y(t+\delta;t+\delta,X(t+\delta;t,\eta;u),u)\big].
\end{equation}

From \eqref{12.13-eq11} and \eqref{12.13-eq12}, we see that
\begin{eqnarray*}
\wt V(t,\eta) \3n& =\3n &
	\operatorname*{ess~inf}_{u(\cdot)\in\mathcal{U}[t,T]}G^{t,\eta;u}_{t,T}\big[\Phi(X(T;t,\eta,u))\big]=
	\operatorname*{ess~inf}_{u(\cdot)\in\mathcal{U}[t,T]}G^{t,\eta;u}_{t,t+\delta}\big[Y(t+\delta;t,x;u)\big]\\
	&  =\3n&
	\operatorname*{ess~inf}_{u(\cdot)\in\mathcal{U}[t,T]}G^{t,\eta;u}_{t,t+\delta}\big[Y(t+\delta;t+\delta,X(t+\delta;t,x;u),u)\big].
\end{eqnarray*}
From the classical comparison theorem of BSDE (e.g., \cite[Theorem 2.2]{KPQ1997}), we have that
$$ 
\wt V(t,\eta)\geq
\operatorname*{ess~inf}_{u(\cdot)\in\mathcal{U}[t,T]}G^{t,\eta;u}_{t,t+\delta}\big[\wt V(t+\delta,X(t+\delta;t,\eta;u))\big].
$$
On the other hand, for every $\varepsilon>0$, we can find an
admissible control $\bar{u}(\cdot)\in\mathcal{U}[t,T]$ such that
$$ 
\wt V(t+\delta,X(t+\delta;t,\eta;\bar{u}))\geq
Y(t+\delta;t+\delta,X(t+\delta;t,\eta,\bar{u}),\bar{u})-\varepsilon.
$$
From this and the comparison theorem of BSDE, we get
$$ 
\wt V(t,\eta)\leq
\operatorname*{ess~inf}_{u(\cdot)\in\mathcal{U}[t,T]}G^{t,\eta;u}_{t,t+\delta}\big[\wt V(t+\delta,X(t+\delta;t,\eta,u))+\varepsilon\big].
$$
From Proposition \ref{prop3.2}, there exists a positive constant $C_0$ such
that
$$ 
\wt V(t,\eta)\leq
\operatorname*{ess~inf}_{u(\cdot)\in\mathcal{U}[t,T]}G^{t,\eta;u}_{t,t+\delta}\big[\wt V(t+\delta,X(t+\delta;t,\eta,u))\big]+C_0\varepsilon.
$$
Therefore, letting $\varepsilon\to 0$, we obtain \eqref{th3.1-eq1}.
\endpf

Theorems \ref{th3.2} and \ref{th3.3} follows from Theorem \ref{th3.1} immediately.

\section{Relationships between PMP and DPP: Smooth Case}\label{sec-smooth}

In this section, we head to the relationship between PMP and DPP when the value function is smooth.

\begin{theorem}\label{th4.1}  
Let {\rm ({\bf S1})--({\bf S3})}  hold and let $ \eta \in H$ be fixed, $(\overline{X} (\cdot), \bar{u}(\cdot), p(\cdot),
q(\cdot))$ be an optimal 4-tuple of \textbf{Problem} $\boldsymbol{(S_{\eta})}$. Suppose that
the value function $V\in C^{1,2}([0,T];H)$ and $V_x$ ranges in $D(A^*)$. Then\vspace{-2mm}
\begin{equation}\label{th4.1-eq1}
\begin{array}{ll}
\ds G\big(t,\overline{X}(t),\bar{u}(t),-V_{x}(t,\overline{X}(t)), -V_{xx}(t,\overline{X}(t))\big)  \\
\ns\ds = \max\limits_{u\in U}G\big(t,\overline{X}(t),u,-V_{x}(t,\overline{X}(t)),-V_{xx}(t,\overline{X}(t))\big), \q \mbox{ a.e. }  (t,\omega)\in [0,T]\times\Omega,
\end{array}
\end{equation}
where $G$ is given in \eqref{prop3.7-eq2}.
Further, if $V\in C^{1,3}([0,T]\times H)$ with $V_{tx}\in H$ and
$A^*V_{xx}\in L(H)$ being continuous on 
$(0,T)\times H$, then\vspace{-2mm}
\begin{equation}\label{th4.1-eq2}
\begin{cases}
V_{x}(t,\overline{X}(t))=-p(t),  \\
V_{xx}(t,\overline{X}(t))b(t,\overline{X}(t),\bar{u}(t))=-q(t),
\end{cases}\q \mbox{ a.e. }  (t,\omega)\in [0,T]\times\Omega.
\end{equation}
\end{theorem}

{\it Proof}. From Theorem \ref{th3.3} and the fact that
$(\overline{X}(\cdot),\bar{u}(\cdot))$ is an optimal pair, we have
$$
V(t,\overline{X}(t))=\mathbb{E}\Big(\int_t^T f(r,\overline{X}(r),\bar{u}(r))dr+h(\overline{X}(T))\Big|\cF_t\Big), \q \forall t\in [0,T],\ \mathbb{P}\mbox{-a.s.}
$$
Let
\begin{equation*}
m(t)\deq\mathbb{E}\Big(\int_0^T f(r,\overline{X}(r),\bar{u}(r))dr+h(\overline{X}(T))\Big|\mathcal{F}_t \Big), \q   t\in [0,T].
\end{equation*}

By  (S4), we have\vspace{-4mm}
\begin{eqnarray*}
  & & \int_0^T \Big[\mathbb{E}\Big|\mathbb{E}\(\int_0^T f(r,\overline{X}(r),\bar{u}(r))dr+h(\overline{X}(T))\Big|\mathcal{F}_t\)\Big|^2\Big]^{1/2}dt\\
  & & \leq T\Big[\mathbb{E}\(\int_0^T |f(r,\overline{X}(r),\bar{u}(r))|^2dr+|h(\overline{X}(T))|^2\)\Big]^{1/2}\\
  & & \leq C\Big\{\mathbb{E}\[\int_0^T \big(1+|\overline{X}(r)|_H^2\big)dr+1+|\overline{X}(T)|_H^2\]\Big\}^{1/2}\\
  & & \leq C\big|\overline{X}(\cdot)\big|_{C_{\mathbb{F}}([0,T];L^2(\Omega;H))},
\end{eqnarray*}
which implies that $m(\cdot) \in L_{\mathbb{F}}^{1}(0,T;L^{2}(\Omega,H))$.
Thus, by the martingale representation theorem (see \cite[Corollary 2.145]{Lu2021} for example), there is a unique $K(\cdot,\cdot,\cdot)\in
L^{1}(0,T;L_{\mathbb{F}}^{2}(0,T;\mathcal{L}_{2}^{0}))$ such that\vspace{-2mm}
$$m(t)=\mathbb{E}m(t)+\int_0^t K(t,r)dW(r), \q t\in [0,T].$$
Then we have\vspace{-4mm}
\begin{eqnarray}\label{12.13-eq49}
   V(t,\overline{X}(t))\3n & =\3n & m(t)-\int_0^t f(r,\overline{X}(r),\bar{u}(r))dr\nonumber\\
   & =\3n &  V(0,\eta)-\int_0^t f(r,\overline{X}(r),\bar{u}(r))dr+\int_0^t K(t,r)dW(r).
\end{eqnarray}
On the other hand, applying It\^o's formula to
$V(t,\overline{X}(t))$, we obtain
\begin{eqnarray}\label{12.13-eq50}
  V(t,\overline{X}(t)) \!\3n& =\3n & V(0,\eta)+ \!\int_0^t\! \Big(V_{r}(r,\overline{X}(r))\!+\! \langle A^*V_x(r,\overline{X}(r)),\overline{X}(r)\rangle_H\! +\!\langle V_{x}(r,\overline{X}(r)), a(r,\overline{X}(r),\bar{u}(r))\rangle_H \nonumber\\
  & & \qq\qq\qq +\frac{1}{2} \langle V_{xx}(r,\overline{X}(r))b (r,\overline{X}(r),\bar{u}(r)),b (r,\overline{X}(r),\bar{u}(r))\rangle_{\mathcal{L}_2^0}\Big)dr  \\
  & & + \int_0^t \lan b (r,\overline{X}(r),\bar{u}(r)) dW(r),V_{x}(r,\overline{X}(r))\ran_H.\nonumber
\end{eqnarray}
From \eqref{12.13-eq49} and  \eqref{12.13-eq50}, we conclude that
\begin{equation}
\begin{cases}\ds
   -f(t,\overline{X}(t),\bar{u}(t))= V_{t}(t,\overline{X}(t))+\langle A^*V_{x}(t,\overline{X}(t)), \overline{X}(t)+a(t,\overline{X}(t),\bar{u}(t))\rangle_H\\
\ns\ds   \ \ \ \   \ \ \ \ \ \ \ \ \ \ \ \ \ \ \ \ \  \ \ \ +\frac{1}{2} \langle V_{xx}b (t,\overline{X}(t),\bar{u}(t)),b (t,\overline{X}(t),\bar{u}(t))\rangle_{\mathcal{L}_2^0},\\
\ns\ds   b (r,\overline{X}(r),\bar{u}(r))V_{x}(r,\overline{X}(r))=K(t,t).
\end{cases}
\end{equation}
This, together with the fact that $V$ is a smooth solution of the
HJB equation, implies \eqref{th4.1-eq1}. In addition,  we have\vspace{-2mm}
\begin{eqnarray*}
  & & G(t,\overline{X}(t),\bar{u}(t),-V_{x}(t,\overline{X}(t)),-V_{xx}(t,\overline{X}(t)))  -\langle A^*V_{x}(t,\overline{X}(t)),\overline{X}(t)\rangle_H-V_{t}(t,\overline{X}(t)) \\
   & & \ge G(t,\eta,\bar{u}(t),-V_{x}(t,\eta),-V_{xx}(t,\eta))-\langle A^{*}V_{x}(t,\eta),\eta\rangle_H-V_{t}(t,\eta), \q \ \forall \eta\in H.
\end{eqnarray*}
Consequently, if $V\in C^{1,3}((0,T);H)$ with $V_{tx}$ and
$A^*V_{xx}$ being also continuous, then\vspace{-2mm}
$$
\frac{\partial}{\partial x}\big(G(t,\eta,\bar{u}(t),-V_{x}(t,\eta),-V_{xx}(t,\eta))-\langle A^{*}V_{x}(t,\eta),\eta\rangle_H-V_{t}(t,\eta)\big)|_{x=\overline{X}(t)}=0.
$$
This implies that\vspace{-2mm}
\begin{eqnarray}
   0\3n & =\3n & V_{tx}(t,\overline{X}(t))+ A^*V_{xx}(t,\overline{X}(t)) \overline{X}(t) +A^{*}V_{x}(t,\overline{X}(t)) \nonumber\\
   & & + a_{x}(t,\overline{X}(t),\bar{u}(t))V_{x}(t,\overline{X}(t))+  V_{xx}(t,\overline{X}(t)) a(t,\overline{X}(t),\bar{u}(t))\nonumber\\
   & & +\frac{1}{2} \sum_{j=1}^\infty V_{xxx}(t,\overline{X}(t))\big( b (t,\overline{X}(t),\bar{u}(t))e_j,b (t,\overline{X}(t),\bar{u}(t))e_j\big) \nonumber\\
   & & +b (t,\overline{X}(t),\bar{u}(t))V_{xx}(t,\overline{X}(t))b_{x} (t,\overline{X}(t),\bar{u}(t))^*+f_{x}(t,\overline{X}(t),\bar{u}(t)),
\end{eqnarray}
where \vspace{-2mm}
$$\frac{\partial}{\partial x}\langle A^{*}V_{x}(t,\eta),\eta\rangle_H=A^*V_{xx}(t,\overline{X}(t)) \overline{X}(t) +A^{*}V_{x}(t,\overline{X}(t))$$
is due to the fact that $A^{*}V_{x}$ is continuous and $A^*$ is
closed, and $\{e_j\}_{j=1}^\infty$ is an orthonormal basis of $\wt H$. One the other hand, by (S4), we have that
\begin{eqnarray}
   dV_{x}(t,\overline{X}(t)) \3n& =\3n & \(V_{tx}(t,\overline{X}(t))+ A^*V_{xx}(t,\overline{X}(t))\overline{X}(t) + V_{xx}(t,\overline{X}(t))a(t,\overline{X}(t),\bar{u}(t)) \nonumber\\
   & & +\frac{1}{2} \sum_{j=1}^\infty V_{xxx}(t,\overline{X}(t))\big( b (t,\overline{X}(t),\bar{u}(t))e_j,b (t,\overline{X}(t),\bar{u}(t))e_j\big)\)dt \nonumber\\
   & & + V_{xx}(t,\overline{X}(t)) b (t,\overline{X}(t),\bar{u}(t))dW(t) \\
   & =\3n & \big(A^{*}V_{x}(t,\overline{X}(t))+ a_{x}(t,\overline{X}(t),\bar{u}(t))V_{x}(t,\overline{X}(t))\nonumber+f_{x}(t,\overline{X}(t),\bar{u}(t)) \\
   & &  +V_{xx}(t,\overline{X}(t))b_{x} (t,\overline{X}(t),\bar{u}(t))^*b (t,\overline{X}(t),\bar{u}(t))\big)dt \nonumber\\
   & & +V_{xx}(t,\overline{X}(t))b (t,\overline{X}(t),\bar{u}(t)) dW(t).\nonumber
\end{eqnarray}
Noting that $-V_{x}(T,\overline{X}(T))=-h_{x}(\overline{X}(T))$, 
 by the uniqueness of the solution to the first-order adjoint
equation, we obtain \eqref{th4.1-eq2}.
\endpf

\begin{remark}
The second equality in \ref{th4.1-eq1} may be regarded as a ``maximum principle" in terms of the derivatives of the value function. It is different from the PMP   in Theorem \ref{maximum p2-1}, where no
derivative of the value function is involved.
\end{remark}

\section{Relationships between PMP and DPP: Nonsmooth Case}\label{sec-nonsmooth}

In Section \ref{sec-smooth}, we dealt with the case when the value function is
smooth enough. However,  in many cases, the
value function associated to a optimal control problem is not
smooth. Hence it is sensible that we drop the smoothness assumptions
of value function in Theorem \ref{th4.1}.   

\subsection{Differential in Spatial Variable}\label{sec-nonsmooth1}

  We shall use the notions of
super and subdifferentials to approach to our desired result. The
main difficulty in carrying out this construction is that in
infinite dimensional setting, the second-order adjoint equation is
an operator-valued BSEE, in this case, the method of duality
relation by means of mild solution breaks down. Thanks to the
introduction of relaxed transposition solution, we are now able to
legitimate our discussion by substituting the duality relation by
the expected relaxed transposition solution. Before stating the main
theory of this section, let us recall the notion of first-order
and second-order super and subdifferentials.

For $v\in C([0,T]\times H)$ and $(t,\eta)\in [0,T)\times H$, the
second-order parabolic superdifferential of $v$ at $(t,\eta)$ is
defined as follows:
\begin{eqnarray*}
   D_{t,x}^{1,2,+}v(t,\eta) \3n& =\3n & \Big\{(r,p,P)\in \mathbb{R}\times H\times \mathcal{S}(H)\Big| \displaystyle\uplim\limits_{\begin{subarray}{1}
   s\downarrow t, s\in [0,T)\\
   y\to \eta
   \end{subarray}}\frac{1}{|s-t|+|\eta-y|_H^2}\\
   & & \Big[v(s,y)-v(t,\eta)-r(s-t)-\langle p,y-\eta\rangle_H-\frac{1}{2}\langle P(y-\eta),y-\eta\rangle_H\Big]\leq 0 \Big\}.
\end{eqnarray*}
Similarly, the second-order parabolic subdifferential of $v$ at
$(t,\eta)$ is defined as follows:
\begin{eqnarray*}
   D_{t,x}^{1,2,-}v(t,\eta) \3n& =\3n & \Big\{(r,p,P)\in \mathbb{R}\times H\times \mathcal{S}(H)\Big| \displaystyle\lowlim\limits_{\begin{subarray}{1}
   s\downarrow t, s\in [0,T)\\
   y\to \eta
   \end{subarray}}\frac{1}{|s-t|+|\eta-y|_H^2}\\
   & & \Big[v(s,y)-v(t,\eta)-r(s-t)-\langle p,y-\eta\rangle_H-\frac{1}{2}\langle P(y-\eta),y-\eta\rangle_H\Big]\ge 0\Big\}.
\end{eqnarray*}
It is important to point out that the limit in $t$ is from the
right. This fits the general irreversibility of evolution equations.

 Let us recall also  the concept of regular conditional probability, which allows us to regard conditional expectations as merely expectations taken with respect to a conditional measure. More details can be found in \cite[Chapter V, Section  8]{Parthasarathy2005}.

\begin{lemma}\label{prop-2.2}
Let $\mathcal{G}$ be a sub-$\sigma$-algebra of $\mathcal{F}$. Then there exists a map $p:  \Omega\times \mathcal{F}\to [0,1]$, called a regular conditional probability given $\mathcal{G}$, such that 

\ss

(i)  for each $\omega\in \Omega$, $p(\omega,\cdot)$ is a probability measure on $\mathcal{F}$;

\ss

(ii)  for each $A\in \mathcal{F}$, the function $p(\cdot,A)$ is $\mathcal{G}$-measurable;

\ss

(iii)  for each $B\in \mathcal{F}$, $p(\omega,B)=\mathbb{P}(B|\mathcal{G})(\omega)=\mathbb{E}(1_B|\mathcal{G})(\omega),\ \mathbb{P}${\rm-a.s.} 

\ss
 
We write $\mathbb{P}(\cdot|\mathcal{G})(\omega)$ for $p(\omega,\cdot)$.
\end{lemma}
\begin{theorem}\label{th5.1}  Suppose {\rm ({\bf S1})--({\bf S3})}  hold. Let $ \eta \in 
H$ be fixed, and $(\overline{X} (\cdot), \bar{u}(\cdot), p(\cdot),
q(\cdot),P(\cdot),
Q^{(\cdot)}, \widehat{Q}^{(\cdot)})$ be an optimal 6-tuple of \textbf{Problem} $\boldsymbol{(S_{\eta})}$. Suppose that
$V\in C([0,T]\times H)$ is the associated value function. Then
\begin{equation}\label{th5.1-eq1}
\{-p(t)\}\times [-P(t),\infty) \subset  D_x^{2,+}V(t,\overline{X}(t)), \q \forall t\in [0,T],\ \mathbb{P}\mbox{-a.s.,} 
\end{equation}
and
\begin{equation}\label{th5.1-eq2}
D_x^{2,-}V(t,\overline{X}(t))\subset  \{-p(t)\}\times (-\infty,-P(t)], \q \forall t\in [0,T],\ \mathbb{P}\mbox{-a.s.}
\end{equation}
\end{theorem}

{\it Proof}. We borrow some ideas in \cite{Cannarsa1992,Cannarsa1996,Cannarsa2015,Yong1999}. The proof is splited into six steps.

{\bf Step  1}.  Fix a $t\in [0,T].$ For any $z\in H$, consider the following SEE:
\begin{equation}\label{12.13-eq56}
\begin{cases}
\ds    dx^z(r)=\big(Ax^z(r)+a(r,x^z(r),\bar{u}(r))\big)dr+b(r,x^z(r),\bar{u}(r))dW(r),\ \ \ \ r\in(t,T],\\
\ns\ds    x^z(t)=z.
\end{cases}
\end{equation}
Set $\xi^z(r)=x^z(r)-\overline{X}(r)$.  In order to reach an estimate for the conditional expectation of $\xi^z(\cdot)$ w.r.t. $\mathcal{F}_t$, we regard equation \eqref{12.13-eq56}  as an SEE
on $(\Omega, \mathcal{F},\textbf{F},
\mathbb{P}(\cdot|\mathcal{F}_t)(\omega))$, for $\mathbb{P}\mbox{-a.s.} \ \omega$. This can be done  by virtue of Proposition \ref{prop-2.2}. 
Similar to the proof of Proposition \ref{prop3.1},  for any $k\ge 1$, we can obtain the following estimate:
\begin{equation}
\mathbb{E}\big(\sup\limits_{t \leq r \leq T} |\xi^z(r)|_H^{2k}\big|
\mathcal{F}_t\big)\leq K|z-\overline{X}(t)|_H^{2k},\q \mathbb{P}\mbox{-a.s.}
\end{equation}
Now we present the equation for $\xi^z(\cdot)$ in two different ways
based on different orders of expansions:
\begin{equation}\label{12.13-eq58}
\!\begin{cases}
\ds  d\xi^z(r)\! =\! \big(A\xi^z(r)\!+\! \bar{a}_{x}(r)\xi^z(r)\big) dr\! + \! \bar{b}_{x}(r)\xi^z(r) dW(r)\! +\!\epsilon_{z,a}(r)dr\!+ \! \epsilon_{z,b}(r)dW(r),
  \; r\in(t,T],\\
\ns\ds  \xi^z(t)=z-\overline{X}(t),
\end{cases}
\end{equation}
and
\begin{equation}\label{12.13-eq59}
\!\begin{cases}\ds
  d\xi^z(r)=\(A\xi^z(r)+ \bar{a}_{x}(r)\xi^z(r)+\frac{1}{2}  \bar{a}_{xx}(r)\big(\xi^z(r),\xi^z(r)\big) \) dr\\
\ns\ds  \ \ \ \ \ \ \ \ \ \  \ + \(  \bar{b}_{x}(r) \xi^z(r)\!+\!\frac{1}{2} \bar{b}_{xx}(r)\big(\xi^z(r),\xi^z(r)\big)\)dW(r)\! + \!\tilde\epsilon_{z,a}dr\!+ \! \tilde\epsilon_{z,b}(r)dW(r),
  \; r\in(t,T],\\
\ns\ds  \xi^z(t)=z-\overline{X}(t),
  \end{cases}
\end{equation}
where for $\f=a,b$,\vspace{-4mm}
\begin{equation}
  \bar{\f}_{x}(r)=\f_{x}(r,\overline{X}(r),\bar{u}(r)),\q 
  \bar{\f}_{xx}(r)=\f_{xx}(r,\overline{X}(r),\bar{u}(r)), 
\end{equation}
and
\begin{equation}
\begin{cases}\ds
  \epsilon_{z,\f}(r)=\int_0^1 \big(\f_{x}(r,\overline{X}(r)+\theta \xi^{z}(r),\bar{u}(r))-\bar{\f}_{x}(r)\big)\xi^{z}(r)d\theta, \\
\ns\ds  \tilde \epsilon_{z,\f}(r)=\int_0^1 (1-\theta)   \big(\f_{xx}(r,\overline{X}(r)+\theta \xi^{z}(r),\bar{u}(r))-\bar{\f}_{xx}(r)\big)\big(\xi^{z}(r),\xi^{z}(r)\big) d\theta.
\end{cases}
\end{equation}

\ss

{\bf Step  2}. This step is devoted to showing that for any $k \ge 1$, there exists
a deterministic continuous and increasing function $\delta :
[0,\infty) \to [0,\infty)$, independent of $z\in H$, with $\delta
(r)=o(r)$, such that
\begin{equation}\label{12.13-eq62}
\mathbb{E}\(\int_t^T |\epsilon_{z,\f}(r)|_H^{2k}dr\Big|\mathcal{F}_t\)(\omega)\leq \delta \big(|z-\overline{X}(t,\omega)|_H^{2k}\big), \hspace{1.06cm} \mathbb{P}\mbox{-a.s.,}
\end{equation}
and
\begin{equation}\label{12.13-eq63}
\mathbb{E}\(\int_t^T |\tilde\epsilon_{z,\f}(r)|_H^{k}dr\Big|\mathcal{F}_t\)(\omega)\leq \delta \big(|z-\overline{X}(t,\omega)|_H^{2k}\big), \hspace{1.06cm} \mathbb{P}\mbox{-a.s.}
\end{equation}

By setting
$\f_{x}(r,\theta)=\f_{x}(r,\overline{X}(r)+\theta \xi^z (r),\bar{u}(r)),$
and using assumption ({\bf S3}), we have
\begin{eqnarray*}
\mathbb{E}\(\int_t^T |\tilde\epsilon_{z,\f}(r)|_H^{k}dr\Big|\mathcal{F}_t\)
\3n  & \leq\3n &  \int_t^T \mathbb{E} \Big(\int_0^1 |\f_{x}(r,\theta)-\bar{\f}_{x}(r)|_{\mathcal{L}(H)}^{2k}d\theta |\xi^z (r)|_H^{2k}\Big|\mathcal{F}_t\Big)dr\\
& \leq\3n & C\int_t^T \mathbb{E} \big(|\xi^z (r)|_H^{4k}\big|\mathcal{F}_t\big)dr\leq C|z-\overline{X}(t,\omega)|_H^{4k}.
\end{eqnarray*}
Thus \eqref{12.13-eq62} follows for $\delta (x)=x^2$. 

Let
$\f_{xx}(r,\theta)=\f_{xx}(r,\overline{X}(r)+\theta \xi^z (r),\bar{u}(r))$.
Similar to the prove of \eqref{12.13-eq62}, we have
\begin{eqnarray*}
\mathbb{E}\(\int_t^T |\tilde\epsilon_{z,\f}(r)|_H^{k}dr\Big|\mathcal{F}_t\)
  &\3n \leq \3n&  \int_t^T \mathbb{E} \(\int_0^1 |\f_{xx}(r,\theta)-\bar{\f}_{xx}(r)|_{\mathcal{L}(H,H;H)}^{k}d\theta |\xi^z (r)|_H^{2k}\Big|\mathcal{F}_t\)dr\\
  &\3n \leq \3n& \int_t^T \big\{\mathbb{E} \big[\bar{\omega}\big(|\xi^z (r)|_H\big)^{2k}\big|\mathcal{F}_t\big]\big\}^{1/2}\big( \mathbb{E} |\xi^z (r)|_H^{4k}\big|\mathcal{F}_t\big)^{1/2}dr\\
  &\3n \leq \3n& C|z-\overline{X}(t,\omega)|_H^{2k}\int_t^T \big\{ \mathbb{E}^{t}[\bar{\omega}(|\xi^z (r)|_H)^{2k}\big|\mathcal{F}_t]\big\}^{1/2}dr.
\end{eqnarray*}
This implies that \eqref{12.13-eq63} holds for some
$\delta (\cdot)$.  Thus, we can pick the largest $\delta(\cdot)$ in the
above calculations, so that \eqref{12.13-eq62}--\eqref{12.13-eq63} follows for all $z\in H$ with
such $\delta(\cdot)$.

\ss

{\bf Step  3}. Let $\bar{f}_{x}(r)=f_{x}(r,\overline{X}(r),\bar{u}(r))$. By It\^o's formula, we get
\begin{eqnarray}
  & & \mathbb{E}\Big(\int_t^T \langle \bar{f}_{x}(r), \xi^z (r)\rangle_H dr+\langle h_{x}(\overline{X}(T)),\xi^z (T)\rangle_H \Big| \mathcal{F}_t \Big) \nonumber\\
  & & =\langle -p(t),\xi^z (t)\rangle_H  -\frac{1}{2}\mathbb{E}\Big[\int_t^T \big(\langle p(r), \xi^z (r)^*\bar{a}_{xx}(r)\xi^z(r) \rangle_H  + \langle q(r), \xi^z (r)^*\bar{b}_{xx}(r)\xi^z(r) \rangle_{\mathcal{L}_2^0} \big)dr\nonumber\\
  & & \ \ \ -\int_t^T \big(\langle p(r),\tilde\epsilon_{z,a}(r)\rangle_H + \langle q(r), \tilde\epsilon_{z,b}(r)\rangle_{\mathcal{L}_2^0} \big)dr\big|\mathcal{F}_t \Big], \ \ \mathbb{P}\mbox{-a.s.}
\end{eqnarray}
 In Definition \ref{def2.1}, relaxed transposition solution is defined by expectation instead of conditional expectation, hence, we use the notion of regular conditional probability and consider the relaxed transposition solution of \eqref{ad-eq2} on probability space $(\Omega, \mathcal{F},\textbf{F},\mathbb{P}(\cdot|\mathcal{F}_t)(\omega))$, for  a.e. $\omega\in\Omega$. For simplicity of notations, we set $\mathbb{E}^{t}_\omega=\mathbb{E}(\cdot |
\mathcal{F}_{t})(\omega)$, which is the probability expectation related to $\mathbb{P}(\cdot|\mathcal{F}_t)(\omega)$, for  a.e. $\omega\in\Omega$. Then we obtain 
\begin{eqnarray}\label{12.13-eq65}
  & & \mathbb{E}^t_\omega\Big(\int_t^T \langle \mathbb{H}_{xx}(r)\xi^z(r),\xi^z(r)\rangle_H dr-\langle h_{xx}(T)\xi^z(T),\xi^z(T)\rangle_H \Big)\nonumber\\
  & & = \langle P(t)\xi^z(t),\xi^z(t)\rangle_H + \mathbb{E}^t_\omega \int_t^T \langle P(r)\epsilon_{z,a}(r),\xi^z(r)\rangle_H dr \nonumber \\
  & &\ \ \  + \mathbb{E}^t_\omega \int_t^T \big(\langle P(r)\xi^z(r),\epsilon_{z,a}(r)\rangle_H + \langle P(r)\bar{b}_{x}(r)\xi^z(r),\epsilon_{z,b}\rangle_{\mathcal{L}_2^0}\big)dr \nonumber \\
  & & \ \ \ +\mathbb{E}^t_\omega \int_t^T \langle P(r)\epsilon_{z,b}(r),\bar{b}_{x}(r)\xi^z(r)+\epsilon_{z,b}(r)\rangle_{\mathcal{L}_2^0}dr \\
  & & \ \ \ +\mathbb{E}^t_\omega\int_t^T \big(\langle \epsilon_{z,b}(r),\hat{Q}^{(t)}(r)\rangle_{\mathcal{L}_2^0}+ \langle Q^{(t)}(r),
  \epsilon_{z,b}(r)\rangle_{\mathcal{L}_2^0}\big)dr,\qquad \mathbb{P}\mbox{-a.s.} \nonumber
\end{eqnarray}

{\bf Step  4}. Compute $V(t,z)-V(t,\overline{X}(t,\omega))$.

Let $M$ be a countable dense subset of $H$. One can
find a subset $\Omega_{0}\subset \Omega$ with
$\mathbb{P}(\Omega_0)=1$ such that for any $\omega_0 \in  \Omega_0$,
\begin{equation*}
\begin{cases}\ds
    V(t,\overline{X}(t,\omega_0))=\mathbb{E}\(\int_t^T f(r,\overline{X}(r),\bar{u}(r))dr+h(\overline{X}(T))\Big|\mathcal{F}_t\)(\omega_0),\\
\ns\ds    \textup{\eqref{12.13-eq56}, \eqref{12.13-eq62}--\eqref{12.13-eq65} hold for any} \ z\in M,\\
\ns\ds      \sup \limits_{s\leq r\leq T}|p(r,\omega_0)|<+\infty,\\
\ns\ds      P(t,\omega_0)\in \mathcal{L}(H),\ P(\cdot,\omega_0)\xi\in L^2(r,T), \ \forall \xi \in L_{\mathcal{F}_r}^2(\Omega;H), \ \forall r\in [t,T].
\end{cases}
\end{equation*}
The first equality holds because of the Bellman's optimal principle,
while the last two conditions are due to the facts that
$$\mathbb{E}\sup\limits_{0\leq r \leq T} |p(r)|_H^2 <+\infty $$
and that $P(\cdot,\cdot)\in
\cP[0,T]$
respectively. Let $\omega_{0}\in \Omega_0$ be fixed, and again set
$\mathbb{E}^{t}_{\omega_0}=\mathbb{E}(\cdot |
\mathcal{F}_{t})(\omega_{0})$. Then for any $z\in M$, by the
definition of value function, we see that
\begin{eqnarray}
  & & V(t,z)-V(t,\overline{X}(t,\omega_{0})) \nonumber\\
  & & \leq \mathbb{E}^{t}_{\omega_0}\Big[\int_t^T \big(f(r,x^{z}(r),\bar{u}(r))-\bar{f}(r)\big)dr+h(x^{z}(T))-h(\overline{X}(T))\Big]\nonumber\\
  & & =\mathbb{E}^{t}_{\omega_0}\Big(\int_t^T \langle \bar{f}_{x}(r), \xi ^{z}(r)\rangle_H dr+\langle h_{x}(\overline{X}(T)),\xi^{z}(T)\rangle_H \Big)  \\
   & & \ \ \  + \frac{1}{2}\mathbb{E}^{t}_{\omega_0} \Big(\int_t^T \langle \bar{f}_{xx}(r)\xi^{z}(r),\xi^{z}(r) \rangle_H dr+\langle h_{xx}(\overline{X}(T))\xi^{z}(T), \xi^{z}(T)\rangle_H \Big) +  o\big(|z-\overline{X}(t,\omega_{0})|_H^2\big). \nonumber 
\end{eqnarray}
From the above discussion, by the definition of the  relaxed transposition solution to \eqref{ad-eq2}, and
the fact that for all $(t,\eta)\in [0,T]\times H$, $V(t,\eta)$  is
deterministic, we have
\begin{eqnarray}\label{12.13-eq67}
  & &\3n\3n\3n V(t,z)-V(t,\overline{X}(t,\omega_0))\nonumber\\
  & &\3n\3n\3n \leq \!-\langle p(t,\omega_0),\xi^z(t,\omega_0)\rangle_H \!-\! \frac{1}{2}\mathbb{E}^t_{\omega_0}\!\Big[\!\int_t^T\!\!\big(\lan p(r),  \bar{a}_{xx}(r)\big(\xi^z\!(r),\xi^z\!(r)\big) \ran_H \! + \!\lan q(r),  \bar{b}_{xx}(r)\big(\xi^z\!(r),\xi^z\!(r)\big) \ran_{\mathcal{L}_2^0} \big)dr\nonumber\\
  & &   -\int_t^T \big(\langle p(r),\tilde \epsilon_{z,a}(r)\rangle_H + \langle q(r), \tilde \epsilon_{z,b}(r)\rangle_{\mathcal{L}_2^0} \big)dr\Big] \nonumber \\
  & &   +\ \frac{1}{2}\mathbb{E}^{t}_{\omega_0} \Big(\int_t^T \langle \bar{f}_{xx}(r)\xi^{z}(r),\xi^{z}(r) \rangle_H dr+\langle h_{xx}(\overline{X}(T))\xi^{z}(T), \xi^{z}(T)\rangle_H \Big) + \ o(|z-\overline{X}(t,\omega_{0})|_H^2)\nonumber\\
  & &\3n\3n\3n  = -\langle p(t,\omega_0),\xi^z(t,\omega_0)\rangle_H -\frac{1}{2}\mathbb{E}^t_{\omega_0}\Big(\int_t^T \langle \mathbb{H}_{xx}(r)\xi^z(r),\xi^z(r)\rangle_H dr - \langle h_{xx}(T)\xi^z(T),\xi^z(T)\rangle_H \Big)\nonumber\\
   & &  -\mathbb{E}^t_{\omega_0} \int_t^T \big(\langle p(r),\tilde \epsilon_{z,a}(r)\rangle_H + \langle q(r), \tilde \epsilon_{z,b}(r)\rangle_{\mathcal{L}_2^0} \big)dr + o(|z-\overline{X}(t,\omega_0)|_H^2)\\
  & & \3n\3n\3n = - \langle p(t,\omega_0),\xi^z(t,\omega_0)\rangle_H - \frac{1}{2}\langle P(t,\omega_0)\xi^z(t,\omega_0),\xi^z(t,\omega_0)\rangle_H\nonumber\\
  & &   -\mathbb{E}^t_{\omega_0} \int_t^T \big(\langle p(r),\tilde \epsilon_{z,a}(r)\rangle_H + \langle q(r), \tilde \epsilon_{z,b}(r)\rangle_{\mathcal{L}_2^0} \big)dr  -\frac{1}{2}\mathbb{E}^t_{\omega_0} \int_t^T \langle P(r)\epsilon_{z,a}(r),\xi^z(r)\rangle_H dr \nonumber\\
  & & -\mathbb{E}^t_{\omega_0} \int_t^T \big(\langle P(r)\xi^z(r),\epsilon_{z,a}(r)\rangle_H + \langle P(r)\bar{b}_{x}(r)\xi^z(r),\epsilon_{z,b}(r)\rangle_{\mathcal{L}_2^0}\big)dr  \nonumber \\
  & &  -\mathbb{E}^t_{\omega_0} \int_t^T \langle P(r)\epsilon_{z,b}(r),\bar{b}_{x}(r)\xi^z(r)+\epsilon_{z,b}(r)\rangle_{\mathcal{L}_2^0}dr \nonumber\\
  & &  -\mathbb{E}^t_{\omega_0} \int_t^T \big(\langle \epsilon_{z,b}(r),\hat{Q}^{(t)}(r)\rangle_{\mathcal{L}_2^0}+ \langle Q^{(t)}(r), \epsilon_{z,b}(r)\rangle_{\mathcal{L}_2^0}\big)dr + \ o(|z-\overline{X}(t,\omega_0)|_H^2).\nonumber
\end{eqnarray}

{\bf Step  5}. In this step, we get rid of terms containing $\epsilon_{z,\cdot}$ and  $\tilde\epsilon_{z,\cdot}$ in \eqref{12.13-eq67}.

By \eqref{12.13-eq62} and \eqref{12.13-eq63}, we have  
\begin{eqnarray}
  & & \Big|\mathbb{E}^t_{\omega_0} \int_t^T \big( \langle p(r),\tilde \epsilon_{z,a}(r)\rangle_H + \langle q(r), \tilde \epsilon_{z,b}(r)\rangle_{\mathcal{L}_2^0} \big)dr\Big|\nonumber\\
  & \leq \3n& \int_t^T \mathbb{E}^t_{\omega_0}\big(|p(r)|_H|\tilde \epsilon_{z,a}(r)|_H+|q(r)|_{\mathcal{L}_2^0}|\tilde \epsilon_{z,b}(r)|_{\mathcal{L}_2^0}\big)dr \\
  & \leq \3n& |p |_{L^2_{\mathbb{F}}(\Omega,C([t,T];H))}\Big(\mathbb{E}^t_{\omega_0}\int_t^T |\tilde \epsilon_{z,a}(r)|_H^2dr\Big)^{1/2} + |q |_{L_{\mathbb{F}}^2(0,T;\mathcal{L}_2^0)} \Big(\mathbb{E}^t_{\omega_0}\int_t^T |\tilde \epsilon_{z,b}(r)|_{\mathcal{L}_2^0}^2dr\Big)^{1/2}\nonumber\\
  & \leq\3n & C\delta(|z-\overline{X}(t,\omega_0)|_H^2)= o(|z-\overline{X}(t,\omega_0)|_H^2).\nonumber
\end{eqnarray}
Noting  that $P$ is a component of the relaxed transposition
solution to the equation \eqref{ad-eq2}, we get that
\begin{eqnarray}
  & &  \Big|\mathbb{E}^t_{\omega_0} \int_t^T \langle P(r)\epsilon_{z,a}(r), \xi^z(r)\rangle_H dr\Big| \nonumber\\
  & &  \leq\Big[\int_t^T \big(\mathbb{E}^t_{\omega_0}|P(r)\epsilon_{z,a}(r)|_{H}^{4/3}\big)^{3/2}dr\Big]^{1/2}\Big(\int_t^T \big(\mathbb{E}^t_{\omega_0}|\xi^z(r)|_{H}^4\big)^{1/2}dr\Big)^{1/2}
  \\
  & & \leq |P |_{\mathcal{L}(L_\mathbb{F}^2(0,T;L^4(\Omega;H));L_{\mathbb{F}}^2(0,T;L^{4/3}(\Omega;H)))}\Big(\mathbb{E}^t_{\omega_0} \int_t^T |\epsilon_{z,a}(r)|_H^4dr \Big)^{1/4} \Big(\sup\limits_{t\leq r\leq T}\mathbb{E}^t_{\omega_0} |\xi^z(r)|_H^4 \Big)^{1/4}\nonumber\\
  & &  =  o(|z-\overline{X}(t,\omega_0)|_H)O(|z-\overline{X}(t,\omega_0 )|_H)= o(|z-\overline{X}(t,\omega_0)|_H^2).\nonumber
\end{eqnarray}
Similarly, we can show that
\begin{eqnarray}
  & &  \Big|\mathbb{E}^t_{\omega_0} \int_t^T \lan P(r)\xi^z(r),\epsilon_{z,a}(r)\ran_H dr \Big|= o\big(|z-\overline{X}(t,\omega_0)|_H^2\big).
\end{eqnarray}
Since $\bar{b}_x$ is bounded, it is also easy to see that
\begin{eqnarray}
  & & \3n\3n  \Big|\mathbb{E}^t_{\omega_0}\int_t^T \langle P(r)\bar{b}_x(r)\xi^z(r),\epsilon_{z,b}(r)\rangle_{\mathcal{L}_2^0} dr\Big|\nonumber\\
  & &\3n\3n \leq \Big[\int_t^T \big(\mathbb{E}^t_{\omega_0}|P(r)\bar{b}_x(r)\xi^z(r)|_{\mathcal{L}_2^0}^{4/3}\big)^{3/2}dr\Big]^{1/2}\Big(\int_t^T \big(\mathbb{E}^t_{\omega_0}|\epsilon_{z,b}(r)|_{\mathcal{L}_2^0}^4\big)^{1/2}dr\Big)^{1/2}\nonumber\\
  & & \3n\3n\leq\! |P |_{\mathcal{L}(L_\mathbb{F}^2(0,T;L^4(\Omega;H));L_{\mathbb{F}}^2(0,T;L^{4/3}(\Omega;H)))} |\bar{b}_x \xi^z |_{L_\mathbb{F}^2(0,T;L^4(\Omega;\mathcal{L}_2^0))}\Big[\!\int_t^T\!\!\big(\mathbb{E}^t_{\omega_0}|\epsilon_{z,b}(r)|_{\mathcal{L}_2^0}^4\big)^{1/2}dr\Big]^{1/2}\nonumber\\
  & &\3n\3n \leq C|\xi^z |_{L_\mathbb{F}^2(0,T;L^4(\Omega;H))}\Big(\int_t^T \big(\mathbb{E}^t_{\omega_0}|\epsilon_{z,b}(r)|_{\mathcal{L}_2^0}^4\big)^{1/2}dr\Big)^{1/2} \\
  & &\3n\3n \leq C(|z-\overline{X}(t,\omega_0)|_H^4)^{1/4}\delta(|z-\overline{X}(t,\omega_0)|_H^4)^{1/4}= o(|z-\overline{X}(t,\omega_0)|_H^2).\nonumber
\end{eqnarray}
Similarly,
\begin{eqnarray}
  & &  \Big|\mathbb{E}^t_{\omega_0} \int_t^T \langle P(r)\epsilon_{z,b}(r),\bar{b}_x(r)\xi^z(r)+\epsilon_{z,b}(r)\rangle_{\mathcal{L}_2^0} dr\Big|=  o(|z-\overline{X}(t,\omega_0)|_H^2).
\end{eqnarray}
Finally, by \eqref{12.13-eq62} and noting that $\hat{Q}^{(t)}$ is part of the
relaxed transposition solution of \eqref{ad-eq2}, we obtain
\begin{eqnarray}
  & &  \Big|\mathbb{E}^t_{\omega_0} \int_t^T \langle \epsilon_{z,b}(r), \hat{Q}^{(t)}(r)\rangle_{\mathcal{L}_2^0}dr\Big|\nonumber\\
   & & \leq \big|\hat{Q}^{(t)}(0,0,\epsilon_{z,b}(\cdot))\big|_{L_\mathbb{F}^2(0,T;L^{4/3}(\Omega;\mathcal{L}_2^0))}\Big[\mathbb{E}^t_{\omega_0}\(\int_t^T |\epsilon_{z,b}(r)|_{\mathcal{L}_2^0}^4dr\)^{1/2}\Big]^{1/2}\nonumber\\
  & & \leq C|\epsilon_{z,b} |_{L_\mathbb{F}^2(0,T;L^4(\Omega,\mathcal{L}_2^0))} \Big[\mathbb{E}^t_{\omega_0}\(\int_t^T |\epsilon_{z,b}(r)|_{\mathcal{L}_2^0}^4dr\)^{1/2}\Big]^{1/2}\\
  & & \leq C \Big[\mathbb{E}^t_{\omega_0}\(\int_t^T |\epsilon_{z,b}(r)|_{\mathcal{L}_2^0}^4dr\)\Big]^{1/2}\leq \big(\delta(|z-\overline{X}(t,\omega_0)|_H^4)\big)^{1/2}\nonumber\\
  & & =o\big(|z-\overline{X}(t,\omega_0)|_H^2\big).\nonumber
\end{eqnarray}
From a similar argument, we arrive at
\begin{eqnarray}
 \Big|\mathbb{E}^t_{\omega_0} \int_t^T \langle Q^{(t)}(r), \epsilon_{z,b}(r)\rangle_{\mathcal{
  L}_2^0}dr\Big|=o\big(|z-\overline{X}(t,\omega_0)|_H^2\big).
\end{eqnarray}

{\bf Step  6}. In this step, we complete the proof.

From the above discussion, we obtain
\begin{eqnarray}\label{12.13-eq75}
  & &\3n\3n\3n V(t,z)-V(t,\overline{X}(t,\omega_0))\nonumber\\
  & &\3n\3n\3n = \langle V_{x}(t,\overline{X}(t,\omega_0)),z-\overline{X}(t,\omega_0) \rangle_H  +\frac{1}{2} \langle V_{xx}(t,\overline{X}(t,\omega_0))(z-\overline{X}(t,\omega_0)),z-\overline{X}(t,\omega_0)\rangle_H \\
  & &\3n\3n\3n \leq - \langle p(t,\omega_0),\xi^z(t,\omega_0)\rangle_H - \frac{1}{2}\langle P(t,\omega_0)\xi^z(t,\omega_0),\xi^z(t,\omega_0)\rangle_H +  o(|z-\overline{X}(t,\omega_0)|_H^2).\nonumber
\end{eqnarray}
Note that the term $o(|z-\overline{X}(t,\omega_0)|_H^2)$ is independent
of $z$. Thus by the continuity of $V(t,\cdot)$, we see that \eqref{12.13-eq75}  holds for all $x\in H$, which proves
$$(-p(t),-P(t))\in D_x^{2,+}V(t,\overline{X}(t)).$$
By the definition of $D_x^{2,+}V(t,\overline{X}(t))$, we obtain \eqref{th5.1-eq1}.

Let us now show \eqref{th5.1-eq2}. Fix an $\omega$ such that \eqref{12.13-eq75} holds for
any $z\in H$. For any $(p,P)\in D_x^{2,-}V(t,\overline{X}(t))$, by the
definition we have
\begin{eqnarray*}
   0 \3n&  \leq \3n & \liminf_{z\to \overline{X}(t)}\frac{V(t,z)-V(t,\overline{X}(t))-\langle p,z-\overline{X}(t)\rangle_H -\frac{1}{2}\langle P(z-\overline{X}(t)),z-\overline{X}(t)\rangle_H}{|z-\overline{X}(t)|^2}\\
   & \leq\3n & \liminf_{z\to \overline{X}(t)}\frac{-\langle p+p(t),z-\overline{X}(t)\rangle_H -\frac{1}{2}\langle (P+P(t))(z-\overline{X}(t)),z-\overline{X}(t)\rangle_H}{|z-\overline{X}(t)|^2},
\end{eqnarray*}
where the last inequality is due to \eqref{12.13-eq75}. Then it is necessary
that
$$p=-p(t),\hspace{1cm} P\leq P(t). $$
This completes the proof.
\endpf

When $V\in C^{1,2}([0,T]\times H),$
\eqref{th5.1-eq1}--\eqref{th5.1-eq2} is reduced to
\begin{equation*}
V_x(t,\overline{X}(t))=-p(t),\q
V_{xx}(t,\overline{X}(t))\leq -P(t).
\end{equation*}

\subsection{Differentials in the Time Variable}\label{sec-nonsmooth2}

In this section, we proceed to studying the super- and
subdifferential of the value function in the time variable $t$ along
an optimal trajectory.  

\begin{theorem}\label{th6.1} 
	Under the assumption of Theorem \ref{th5.1}, for any
$t$ such that $\overline{X}(t)\in D(A)$ or $p(t)\in D(A^*)$, we have
$$\langle\langle  A\overline{X}(t),p(t)\rangle\rangle + \mathcal{H}(t,\overline{X}(t),\bar{u}(t))\in D_{t+}^{1,+}V(t,\overline{X}(t)), \ \ a.e. \ t\in [0,T],\ \ \mathbb{P}\mbox{-a.s.}$$
where \vspace{-2mm}
$$\langle\langle  A\overline{X}(t),p(t)\rangle\rangle = \begin{cases}\langle  A\overline{X}(t),p(t)\rangle_H, & \mbox{ if }\;\overline{X}(t)\in D(A),\\
	\ns\ds \langle \overline{X}(t), A^*p(t)\rangle_H, & \mbox{ if }\;p(t)\in D(A^*),
	\end{cases}
$$
and
$\mathcal{H}(t,\eta,u)=G(t,\eta,u,p(t),P(t))+\langle b
(t,\eta,u),q(t)-P(t)b (t,\eta,\bar{u}(t))\rangle_{\cL_2^0}$.
\end{theorem}

{\it Proof}. For any $t\in (0,T)$, take $\tau \in (t,T]$. Denote by
$x_{\tau}$ the solution to the following SEE:
\begin{equation}
\begin{cases}
\ds   dx_{\tau}(r)=Ax_{\tau}(r)dr+a(r,x_{\tau}(r),\bar{u}(r))dr+b(r,x_{\tau}(r),\bar{u}(r))dW(r),\ \ r\in (\tau, T],\\
\ns\ds   x_{\tau}(\tau)=\overline{X}(t).
\end{cases}
\end{equation}
Set $\xi_{\tau}(r)=x_{\tau}(r)-\overline{X}(r)$ for $r\in[\tau,T]$.
Working under $\mathbb{P}(\cdot|\mathcal{F}_{\tau})(\omega),\ \mathbb{P}\mbox{-a.s.}\ \omega $, we have
\begin{equation}
    \mathbb{E}\(\sup\limits_{\tau \leq r \leq T} |\xi_{\tau}(r)|_H^{2k}| \mathcal{F}_{\tau}\)\leq C\big|\overline{X}(\tau)-\overline{X}(t)\big|_H^{2k} ,\ \ \ \mathbb{P}\mbox{-a.s.}
\end{equation}
Taking $\mathbb{E}(\cdot|\mathcal{F}_t)$ on both sides and noting
that $\mathcal{F}_t\subset \mathcal{F}_{\tau}$, we obtain
\begin{equation}\label{12.13-eq6.25}
    \mathbb{E}\(\sup\limits_{\tau \leq r \leq T} |\xi_{\tau}(r)|_H^{2k}| \mathcal{F}_t\)\leq C|\tau-t|^k\big(|A\overline{X}(t)|_H+1+|\overline{X}(t)|_H\big)^k\leq C|\tau-t|^k,\ \ \ \mathbb{P}\mbox{-a.s.}
\end{equation}
From the definition of $\xi_{\tau}(\cdot)$, we know that it  satisfies the following variational
equations:
\begin{eqnarray}\label{12.14-eq1}
\3n\!\!\begin{cases}
  d\xi_{\tau}(r)\! =\! \big(A\xi_{\tau}(r)\!+\! \bar{a}_{x}(r)\xi_{\tau}(r)\big) dr\! +\!  \bar{b}_{x}(r) \xi_{\tau}(r) dW\!(r) \!+\!\epsilon_{\tau,a}(r)dr\!+ \! \epsilon_{\tau,b}(r)dW\!(r),
  \,\; r\in(\tau,T],\\
\ns\ds  \xi_{\tau}(\tau)=-[S(\tau-t)-I]\overline{X}(t)-\int_t^{\tau} S(\tau-r)\bar{a}(r)dr-\int_t^{\tau} S(\tau-r)\bar{b}(r)dW(r),
\end{cases}
\end{eqnarray}
and
\begin{eqnarray}\label{12.14-eq2}
\3n\!\begin{cases}\ds
  d\xi_{\tau}(r)=\(A\xi_{\tau}(r)+ \bar{a}_{x}(r)\xi_{\tau}(r)+\frac{1}{2}  \bar{a}_{xx}(r)\big(\xi_{\tau}(r),\xi_{\tau}(r) \big)  \) dr\\
 \ns\ds \ \ \ \ \ \ \ \ \ \ \   +  \(  \bar{b}_{x}(r) \xi_{\tau}(r)\!+\!\frac{1}{2} \bar{b}_{xx}(r)\big(\xi_{\tau}(r),\xi_{\tau}(r) \big)\)dW\!(r)\! +\!   \tilde\epsilon_{\tau,a}(r)d\!+\!  \tilde\epsilon_{\tau,b}(r)dW\!(r), 
\,\;r\in(\tau,T],\\
\ns\ds  \xi_{\tau}(\tau)=-[S(\tau-t)-I]\overline{X}(t)-\int_t^{\tau} S(\tau-r)\bar{a}(r)dr-\int_t^{\tau} S(\tau-r)\bar{b}(r)dW(r).
  \end{cases}
\end{eqnarray}
Here for $\f=a,b$,
\begin{equation}\label{12.13-eq6.29}
	\begin{cases}\ds
		\epsilon_{\tau,\f}(r)=\int_0^1 \big(\f_{x}(r,\overline{X}(r)+\theta \xi_{\tau}(r),\bar{u}(r))-\bar{\f}_{x}(r)\big)\xi_{\tau}(r)d\theta, \\
		\ns\ds  \tilde \epsilon_{\tau,\f}(r)=\int_0^1 (1-\theta)\xi_{\tau}(r)^*  \big(\f_{xx}(r,\overline{X}(r)+\theta \xi_{\tau}(r),\bar{u}(r))-\bar{\f}_{xx}(r)\big)\xi_{\tau}(r) d\theta.
	\end{cases}
\end{equation}
Similar to \eqref{12.13-eq62} and \eqref{12.13-eq63}, one can prove that for any $k\ge 1$,
\begin{equation}
\begin{cases}\ds
  \mathbb{E}\(\int_{\tau}^T |\epsilon_{\tau,\f}(r)|_H^{2k}dr|\mathcal{F}_t \) \leq \delta (|\tau-t|^k), \qquad \mathbb{P}\mbox{-a.s.},\\
\ns\ds    \mathbb{E}\(\int_{\tau}^T |\tilde \epsilon_{\tau,\f}(r)|_H^{k}dr|\mathcal{F}_t \) \leq \delta (|\tau-t|^k), \qquad \mathbb{P}\mbox{-a.s.}, 
\end{cases}
\end{equation}
for some deterministic continuous function $\delta:[0,\infty)\to
[0,\infty)$ with $\frac{\delta(r)}{r}\to 0$ as $r\to 0$. By the definition of the value function $V$, we have
\begin{equation}\label{12.13-eq83}
    V(\tau,\overline{X}(t))\leq \mathbb{E}\Big(\int_{\tau}^T f(r,x_{\tau}(r),\bar{u}(r))dr+h(x_{\tau}(T))\Big|\mathcal{F}_{\tau} \Big),\ \ \ \mathbb{P}\mbox{-a.s.}
\end{equation}
Taking $\mE(\cdot|\mathcal{F}_t^s)$ on both sides of \eqref{12.13-eq83} and noting
that $t\leq \tau$, we conclude that
\begin{equation}\label{12.13-eq6.31}
    V(\tau,\overline{X}(t))\leq \mathbb{E}\Big(\int_{\tau}^T f(r,x_{\tau}(r),\bar{u}(r))dr+h(x_{\tau}(T))\Big|\mathcal{F}_t \Big),\ \ \ \mathbb{P}\mbox{-a.s.}
\end{equation}
Choose a subset $\Omega_0\in \Omega$, with $\mathbb{P}(\Omega_0)=1$
such that for any $\omega_0\in \Omega_0$,
\begin{equation}
    \begin{cases}\ds
      V(t,\overline{X}(t,\omega_0))=\mathbb{E}\Big(\int_t^T f(r,\overline{X}(r),\bar{u}(r))dr+h(\overline{X}(T))\Big|\mathcal{F}_t \Big)(\omega_0),\nonumber\\
      \textup{\eqref{12.13-eq6.25}, \eqref{12.13-eq6.29} and \eqref{12.13-eq6.31} are satisfied for any rational}\  \tau>t,\nonumber\\
      \sup \limits_{s\leq r\leq T}|p(r,\omega_0)|<+\infty,\\
    P(t,\omega_0)\in L(H),\ P(\cdot,\omega_0)\xi\in L^2(r,T;H), \ \forall \xi \in L_{\mathcal{F}_r}^2(\Omega;H), \ \forall r\in [0,T].
    \end{cases}
\end{equation}
Let $\omega_0\in \Omega_0$ be fixed, and set
$\mathbb{E}^t_{\omega_0}=\mathbb{E}(\cdot|\mathcal{F}_t ) (\omega_0)$. Then
for any rational $\tau>t,$ we have
\begin{eqnarray}
   & & V(\tau,\overline{X}(t,\omega_0))-V(t,\overline{X}(t,\omega_0))\nonumber\\
   & \leq\3n & \mathbb{E}^t_{\omega_0}\Big\{-\int_t^{\tau}\bar{f}(r)dr+\int_{\tau}^T\big[f(r,x_{\tau}(r),\bar{u}(r))-\bar{f}(r)\big]dr+h(x_{\tau}(T))-h(\overline{X}(T))\Big\}\nonumber\\
   & =\3n & \mathbb{E}^t_{\omega_0}\Big(-\int_t^{\tau}\bar{f}(r)dr+\int_{\tau}^T \langle \bar{f}_x(r),\xi_{\tau}(r)\rangle_H dr+\langle h_x(\overline{X}(T)),\xi_{\tau}(T)\rangle_H  \\
   & & \qq+\frac{1}{2}\int_{\tau}^T \langle \bar{f}_{xx}(r)\xi_{\tau}(r),\xi_{\tau}(r)\rangle_H dr+\frac{1}{2}\langle h_{xx}(\overline{X}(T))\xi_{\tau}(T),\xi_{\tau}(T)\rangle_H \Big)+o(|\tau-t|)\nonumber.
\end{eqnarray}
Similar to the argument in the forth and fifth steps in the proof of Theorm
\ref{th6.1}, we can obtain
\begin{eqnarray}\label{12.13-eq86}
   & & V(\tau,\overline{X}(t,\omega_0))-V(t,\overline{X}(t,\omega_0))\nonumber\\
   & \leq\3n   & -\mathbb{E}^t_{\omega_0}\int_t^{\tau}\bar{f}(r)dr-\mathbb{E}^t_{\omega_0}\(\langle p(\tau),\xi_{\tau}(\tau)\rangle_H - \frac{1}{2}\langle P(\tau)\xi_{\tau}(\tau),\xi_{\tau}(\tau)\rangle_H\)\nonumber\\
  & &     -\mathbb{E}^t_{\omega_0} \int_{\tau}^T \big(\langle p(r),\tilde\epsilon_{\tau,a}(r)\rangle_H + \langle q(r), \tilde\epsilon_{\tau,b}(r)\rangle_{\mathcal{L}_2^0} \big)dr  -\frac{1}{2}\mathbb{E}^t_{\omega_0} \int_{\tau}^T \langle P(r)\epsilon_{\tau,a}(r),\xi_{\tau}(r)\rangle_H dr  \nonumber\\
  & & -\mathbb{E}^t_{\omega_0} \int_{\tau}^T \big(\langle P(r)\xi_{\tau}(r),\epsilon_{\tau,a}(r)\rangle_H + \langle P(r)\bar{b}_{x}(r)\xi_{\tau}(r),\epsilon_{\tau,b}(r)\rangle_{\mathcal{L}_2^0}\big)dr\\
  & & -\mathbb{E}^t_{\omega_0} \int_{\tau}^T \langle P(r)\epsilon_{\tau,b}(r),\bar{b}_{x}(r)\xi_{\tau}(r)+\epsilon_{\tau,b}(r)\rangle_{\mathcal{L}_2^0}dr  \nonumber \\
  & & -\mathbb{E}^t_{\omega_0} \int_{\tau}^T \big(\langle \epsilon_{\tau,b}(r),\hat{Q}^{(\tau)}(r)\rangle_{\mathcal{L}_2^0}+ \langle Q^{(\tau)}(r), \epsilon_{\tau,b}(r)\rangle_{\mathcal{L}_2^0}\big)dr + \ o(|\tau-t|)\nonumber\\
  & =\3n   & -\mathbb{E}^t_{\omega_0}\int_t^{\tau}\bar{f}(r)dr-\mathbb{E}^t_{\omega_0}\( \langle p(\tau),\xi_{\tau}(\tau)\rangle_H - \frac{1}{2}\langle P(\tau)\xi_{\tau}(\tau),\xi_{\tau}(\tau)\rangle_H\) +o(|\tau-t|).\nonumber
\end{eqnarray}
Now, let us estimate the terms on the right-hand side of \eqref{12.13-eq86}. To
this end, we first note that for any $\phi,\phi'\in
L_{\mathbb{F}}^2(0,T;H)$, and $\psi\in
L_{\mathbb{F}}^2(0,T;\mathcal{L}_2^0) $, it holds that
\begin{eqnarray}\label{12.13-eq87}
    \mathbb{E}^t_{\omega_0}\Big\langle \int_t^{\tau} \phi(r)dr,\int_t^{\tau} \phi'(r)dr\Big\rangle_H\3n &\leq\3n &  \Big(\mathbb{E}^t_{\omega_0}\Big|\int_t^{\tau} \phi(r)dr\Big|_H^2\Big)^{\frac{1}{2}} \Big(\mathbb{E}^t_{\omega_0}\Big| \int_t^{\tau} \phi'(r)dr\Big|_H^2\Big)^{\frac{1}{2}}\nonumber\\
   & =\3n &(\tau-t)\Big(\int_t^{\tau}\mathbb{E}^t_{\omega_0}|\phi(r)|_H^2dr\int_t^{\tau}\mathbb{E}^t_{\omega_0}|\phi'(r)|_H^2dr\Big)^{\frac{1}{2}}\\
   & =\3n & o(|\tau-t|), \hspace{1cm}\textup{as}\ \tau\downarrow t,\ \ \forall \ t\in [0,T),\ \ \mathbb{P}\mbox{-a.s}.\nonumber
\end{eqnarray}
and due to the full Lebesgue measure for integrable functions and
the fact that $t\mapsto \mathcal{F}_t$ is continuous w.r.t $t$, we
have
\begin{eqnarray}\label{12.13-eq88}
     \mathbb{E}^t_{\omega_0}\Big\langle\! \int_t^{\tau}\! \phi(r)dr,\int_t^{\tau}\! \psi(r)dW(r)\Big\rangle_H \3n&\leq\3n&  \Big(\mathbb{E}^t_{\omega_0}\Big|\int_t^{\tau} \phi(r)dr\Big|_H^2\Big)^{\frac{1}{2}} \Big(\mathbb{E}^t_{\omega_0}\Big|\int_t^{\tau} \psi(r)dW(r)\Big|_{\mathcal{L}_2^0}^2\Big)^{\frac{1}{2}}\nonumber\\
   & =\3n &(\tau-t)^{\frac{1}{2}}\Big(\!\int_t^{\tau}\!\mathbb{E}^t_{\omega_0}|\phi(r)|_H^2dr\!\int_t^{\tau}\!\mathbb{E}^t_{\omega_0}|\psi(r)|_{\mathcal{L}_2^0}^2dr\Big)^{\frac{1}{2}}\\
   & =\3n & o(|\tau-t|), \hspace{1cm}\textup{as}\ \tau\downarrow t,\ \ a.e.\  t\in [0,T),\ \ \mathbb{P}\mbox{-a.s}.\nonumber
\end{eqnarray}
Thus, by \eqref{12.13-eq87} and \eqref{12.13-eq88},
\begin{eqnarray}\label{12.13-eq88-1}
  & &\3n\3n \3n\3n\3n\3n   \mathbb{E}^t_{\omega_0} \langle p(\tau),\xi_{\tau}(\tau) \rangle_H =\mathbb{E}^t_{\omega_0}\big(\langle p(t),\xi_{\tau}(\tau)\rangle_H+\langle p(\tau)-p(t),\xi_{\tau}(\tau)\rangle_H\big)\nonumber\\
   &  &\3n\3n \3n\3n\3n\3n   =  \mathbb{E}^t_{\omega_0}\Big\{\Big\langle p(t), -[S(\tau-t)-I]\overline{X}(t)-\int_t^{\tau} S(\tau-r)\bar{a}(r)dr-\int_t^{\tau} S(\tau-r)\bar{b}(r)dW(r)\Big\rangle_H\nonumber\\
   & & +\Big\langle [S(\tau-t)-I]S(T-\tau)h_{x}(\overline{X}(T))-\int_t^{\tau}S(r-t)\big[\bar{a}_x(r)^*p(r)+\bar{b}_x(r)^*q(r)-\bar{f}_x(r)\big]dr\nonumber\\
   & &   -\int_{\tau}^T[S(\tau-t)-I]S(r-\tau)\big[\bar{a}_x(r)^*p(r) +\bar{b}_x(r)^*q(r) -\bar{f}_x(r)\big]dr\\
   & &  
   + \int_t^{\tau}S(r-t)q(r)dW(r)-\int_{\tau}^T[S(\tau-t)-I]S(r-\tau)q(r)dW(r),\nonumber\\
   & &   -[S(\tau-t)-I]\overline{X}(t)-\int_t^{\tau} S(\tau-r)\bar{a}(r)dr
   -\int_t^{\tau} S(\tau-r)\bar{b}(r)dW(r)\Big\rangle_H\Big\}\nonumber\\
   &  &\3n\3n \3n\3n \3n\3n  =\mathbb{E}^t_{\omega_0}\[-\Big\langle A^*p(t),(\tau-t) \overline{X}(t)- \Big\langle p(t),\int_t^{\tau} S(\tau-r)\bar{a}(r)dr\Big\rangle_H \nonumber\\
   &  & -\int_t^{\tau} \lan S(r-t)q(r), S(\tau-r)\bar{b}(r)\ran_{\cL_2^0}dr\] +o(|\tau-t|).\nonumber 
\end{eqnarray}
By the definition of $\xi_{\tau}(\tau)$, we have
\begin{eqnarray}\label{12.13-eq90}
   & & \mathbb{E}^t_{\omega_0}\langle P(\tau)\xi_{\tau}(\tau),\xi_{\tau}(\tau)\rangle_H\nonumber\\
   & =\3n &  \mathbb{E}^t_{\omega_0}\Big\langle P(\tau)\Big\{[S(\tau-r)-I]\overline{X}(t)-\int_t^{\tau}S(\tau-r)\bar{a}(r)dr-\int_t^{\tau}S(\tau-r)\bar{b}(r)dW(r)\Big\},\nonumber\\
   & & \q\; [S(\tau-r)-I]\overline{X}(t)-\int_t^{\tau}S(\tau-r)\bar{a}(r)dr-\int_t^{\tau}S(\tau-r)\bar{b}(r)dW(r)\Big\rangle_H\nonumber\\
   & =\3n & \mathbb{E}^t_{\omega_0}\int_t^{\tau} \langle P(\tau)\bar{b}(r),\bar{b}(r)\rangle_{\mathcal{L}_2^0} dr+o(|\tau-t|)\\
   & =\3n & \mathbb{E}^t_{\omega_0}\int_t^{\tau}  \langle P(\tau)\bar{b}(r),\bar{b}(r)-\bar{b}(t)\rangle_{\mathcal{L}_2^0} dr+\mathbb{E}^t_{\omega_0}\int_t^{\tau} \big\langle P(\tau)\big(\bar{b}(r)-\bar{b}(t)\big),\bar{b}(t)\big\rangle_{\mathcal{L}_2^0} dr\nonumber\\
   & & +  \mathbb{E}^t_{\omega_0}\int_t^{\tau} \big\langle \big(P(\tau)-P(t)\big)\bar{b}(t),\bar{b}(t)\big\rangle_{\mathcal{L}_2^0} dr+\mathbb{E}^t_{\omega_0}\int_t^{\tau} \langle P(t)\bar{b}(t),\bar{b}(t)\rangle_{\mathcal{L}_2^0} dr +o(|\tau-t|).\nonumber
\end{eqnarray}
Let us estimate the first three terms on the right-hand side of
\eqref{12.13-eq90}. Since $P(\cdot)\bar{b}(r)\in
D_{\mathbb{F}}([r,T];$ $L^{4/3}(\Omega;\mathcal{L}_2^0))$, we see that
\begin{eqnarray}\label{12.13-eq91}
   & & \mathbb{E}^t_{\omega_0}\int_t^{\tau}  \langle P(\tau)\bar{b}(r),\bar{b}(r)-\bar{b}(t)\rangle_{\mathcal{L}_2^0} dr\nonumber\\
   & \leq\3n & |P(\cdot)\bar{b}(\cdot)|_{D_{\mathbb{F}}([t,T];L^{4/3}(\Omega;\mathcal{L}_2^0))}\int_t^{\tau}\big(\mathbb{E}^t_{\omega_0}|\bar{b}(r)-\bar{b}(t)|_{\mathcal{L}_2^0}^4\big)^{\frac{1}{4}}dr\\
   & =\3n & o(|\tau-t|), \hspace{1cm}  \textup{as}\ \tau\downarrow t,\ \ a.e.\  t\in [0,T).\nonumber
\end{eqnarray}
Similarly,
\begin{eqnarray}\label{12.13-eq92}
   & & \mathbb{E}^t_{\omega_0}\int_t^{\tau} \big\langle P(\tau)\big(\bar{b}(r)-\bar{b}(t)\big),\bar{b}(t)\big\rangle_{\mathcal{L}_2^0} dr\nonumber\\
   & \leq\3n & \int_t^\tau \big[\mathbb{E}^t_{\omega_0}\big|P(\tau)\big(\bar{b}(r)-\bar{b}(t)\big)\big|_{\mathcal{L}_2^0}^{\frac{4}{3}}\big]^{\frac{3}{4}}dr|\bar{b}(t)|_{L^4_{\mathcal{F}_t}(\Omega;\mathcal{L}_2^0)}\\
   & =\3n & o(|\tau-t|), \hspace{1cm}  \textup{as}\ \tau\downarrow t,\ \ a.e.\  t\in [0,T).\nonumber
\end{eqnarray}
and
\begin{eqnarray}\label{12.13-eq93}
   & & \mathbb{E}^t_{\omega_0}\int_t^{\tau} \big\langle \big(P(\tau)-P(t)\big)\bar{b}(t),\bar{b}(t)\big\rangle_{\mathcal{L}_2^0} dr\nonumber\\
   & =\3n & (\tau-t)\mathbb{E}^t_{\omega_0} \big\langle \big(P(\tau)-P(t)\big)\bar{b}(t),\bar{b}(t)\big\rangle_{\mathcal{L}_2^0} \nonumber\\
   & \leq\3n &  (\tau-t)\big[\mathbb{E}^t_{\omega_0}\big|\big(P(\tau)-P(t)\big)\bar{b}(t)\big|^{\frac{4}{3}}_{\mathcal{L}_2^0}\big]^{\frac{3}{4}}\big|\bar{b}(t)\big|_{L_{\mathcal{F}_t}^4(\Omega;\mathcal{L}_2^0)}\\
   & \leq\3n & (\tau-t)\big|(P(\tau)-P(t)\big)\bar{b}(t)\big|_{L^{4/3}_{\mathcal{F}_t}(\Omega;\mathcal{L}_2^0)}\big|\bar{b}(t)|_{L^4_{\mathcal{F}_t}(\Omega;\mathcal{L}_2^0)}\nonumber\\
   & =\3n & o(|\tau-t|), \q\textup{as}\ \tau\downarrow t, \nonumber
\end{eqnarray}
where the last equality if due to the right continuity of $P(\cdot)\xi$ in
$L^{4/3}_{\mathcal{F}_t}(\Omega;H)$, $\forall \xi \in L_{\mathcal{F}_t}^4(\Omega;H)
$. Combining \eqref{12.13-eq90}--\eqref{12.13-eq93}, we see that
\begin{eqnarray}\label{12.13-eq93-1}
   &&\mathbb{E}^t_{\omega_0}\langle P(\tau)\xi_{\tau}(\tau),\xi_{\tau}(\tau)\rangle_H 
 \nonumber \\
  & =\3n & \mathbb{E}^t_{\omega_0}\int_t^{\tau} \langle P(t)\bar{b}(t),\bar{b}(t)\rangle_{\mathcal{L}_2^0} dr +o(|\tau-t|)\\
   & = \3n& (\tau\!-t)\langle P(t)\bar{b}(t),\bar{b}(t)\rangle_{\mathcal{L}_2^0} +o(|\tau-t|),  \q\textup{as}\ \tau\downarrow t,\ \ a.e.\  t\in [0,T).\nonumber
\end{eqnarray}
It follows from \eqref{12.13-eq86}, \eqref{12.13-eq88-1} and  \eqref{12.13-eq93-1} that for any rational $\tau>t$
and at $\omega=\omega_0$,
\begin{eqnarray*}
   & &\3n\3n\3n\3n\3n\3n V(\tau,\overline{X}(t))-V(t,\overline{X}(t))\nonumber\\
   &&\3n\3n\3n\3n\3n\3n \leq \mathbb{E}^t_{\omega_0}\Big\{\big\langle p(t),\big[S(\tau-t)-I\big]\overline{X}(t)\big\rangle_H+\Big\langle p(t),\int_t^{\tau}S(\tau-r)\bar{a}(r)dr\Big\rangle_H\nonumber\\
   & & +\int_t^{\tau}\! \lan S(r\!-\!t)q(r), S(\tau\!-\!r)\bar{b}(r)\ran_{\cL_2^0}dr-\frac{1}{2}(\tau\!-t)\langle P(t)\bar{b}(t),\bar{b}(t)\rangle_{\mathcal{L}_2^0} -\int_t^{\tau}\!\bar{f}(r)dr\Big\} +o(|\tau\!-t|)\nonumber\\
   &&\3n\3n\3n\3n\3n\3n \leq  (\tau-t)\big[\lan\lan A\overline{X}(t),p(r)\ran\ran+\mathcal{H}(t,\overline{X}(t),\bar{u}(t))\big]+o(|\tau-t|),
\end{eqnarray*}
which completes the proof.

\subsection{Examples}\label{sec-exam}

In this section, we present two examples which fulfills the assumptions in Theorem \ref{th5.1} and/or \ref{th6.1}.

Let $G\subset\dbR^n$ be a bounded domain with the smooth boundary $\pa G$. Let $H=L^2(G)$ and $U$  be a bounded close subset of $L^2(G)$.   
Consider the following
stochastic parabolic equation:
\begin{equation}\label{system3}
\begin{cases}
\ds   dy =\big(\Delta y + \tilde a(t,y,u)\big) dt+\tilde b(t,y,u) dW(t) &\textup{in }  (0,T]\times G,\\
\ns\ds   y =0 &\textup{on }   (0,T]\times \pa G,\\
\ns\ds   y(0)=\eta &\textup{in }  G.
\end{cases}
\end{equation}
where $ \eta \in L^2(G)$, $u(\cdot)\in \cU[0,T]$, and $\tilde a$ and $\tilde b$ satisfy the following condition:

\ss

\no{\bf (B1)} {\it For $\f=\tilde a,\tilde b$, suppose that
	$\f(\cd,\cd,\cd):[0,T]\times \dbR\times \dbR\to
	\dbR$ satisfies : i) For any $(r,u)\in
	\dbR\times \dbR$, the functions
	$\f(\cd,r,u):[0,T] \to \dbR$ is Lebesgue
	measurable; ii) For any $(t,r)\in [0,T]\times 
	\dbR$, the functions $\f(t,r,\cd):\dbR\to \dbR$
	is continuous; and iii) For all
	$(t,r_1,r_2,u)\in [0,T]\times \dbR\times
	\dbR\times \dbR$,
	\begin{equation}\label{ab0}
		\left\{
		\begin{array}{ll}\ds
			|\f(t,r_1,u) - \f(t,r_2,u)|   \leq
			\cC |r_1-r_2|,\\
			\ns\ds |\f(t,0,u)| \leq \cC;
		\end{array}
		\right.
\end{equation} 
 iv)  For all
 $(t,u)\in [0,T]\times  \dbR$,
	$\f(t,\cd,u)$ are $C^2$, and for any $(r,u)\in
	\dbR\times \dbR$ and a.e. $t\in [0,T]$,
	\begin{equation*}\label{ab1}
			|\f_r(t,r,u)| + |\f_{rr}(t,r,u)|
		 \leq \cC.
\end{equation*}}

 Consider the following 
cost functional:
\begin{equation}\label{cost3}
\mathcal{J}(\eta;u(\cdot))= \mathbb{E}\Big[\int_0^T\int_G \tilde f(t,y(t),u(t))dxdt+\int_G \tilde h(y(T))dx 
\Big],
\end{equation}
where $\tilde f$ and $\tilde h$ satisfy the following condition:

\ss

\no{\bf (B2)} {\it    
	$\tilde f(t,\cd,u)$ and $\tilde h(\cd)$ are $C^2$, such that 
	$\tilde f_r(t,r,\cd)$ and $\tilde f_{rr}(t,r,\cd)$ are continuous, and for any $(r,u)\in
	\dbR\times \dbR$ and a.e. $t\in [0,T]$,
	\begin{equation*}\label{ab1-1}
			|\tilde f_r(t,r,u)| +|\tilde h_r(r) |+|\tilde f_{rr}(t,r,u)|
			+|\tilde h_{rr}(r)| \leq \cC.
\end{equation*}}

Under {\bf (B1)} and {\bf (B2)}, it is easy to see that ({\bf S1})--({\bf S3}) hold. Then we know that all assumptions in Theorem \ref{th5.1} are fulfilled. By the regularity of backward stochastic parabolic equations (e.g.,\cite{Du2012}), we know that $A^*p(t)\in L^2(G)$ for a.e. $(t,\om)\in (0,T)\times \Om$. Hence,  all assumptions in Theorem \ref{th6.1} are fulfilled.

\ss

Next, let $H=H_0^1(G)\times L^2(G)$ and $U$  be a bounded close subset of $L^2(G)$.   
Consider the following
stochastic hyperbolic equation
\begin{equation}\label{system4}
	\begin{cases}
		\ds   dy_t =\big(\Delta y + \tilde a(t,y,u)\big) dt+\tilde b(t,y,u) dW(t) &\textup{in }  (0,T]\times G,\\
		\ns\ds   y =0 &\textup{on }   (0,T]\times \pa G,\\
		\ns\ds   y(0)=\eta_1,\; y_t(0)=\eta_2 &\textup{in }  G.
	\end{cases}
\end{equation}
where $ (\eta_1,\eta_2)\in H_0^1(G)\times L^2(G)$, $u(\cdot)\in \cU[0,T]$, and $\tilde a$ and $\tilde b$ satisfy {\bf (B1)}.

 Consider the following 
cost functional:\vspace{-2mm}
\begin{equation}\label{cost4}
	\mathcal{J}(\eta_1,\eta_2;u(\cdot))= \mathbb{E}\Big[\int_0^T\int_G \tilde f(t,y(t),u(t))dxdt+\int_G \tilde h(y(T))dx 
	\Big],
\end{equation}
where $\tilde f$ and $\tilde h$ satisfy  {\bf (B2)}. Under {\bf (B1)} and {\bf (B2)}, it is easy to see that ({\bf S1})--({\bf S3}) hold. Then we know that all assumptions in Theorem \ref{th5.1} are fulfilled. 

\appendix

\end{document}